\documentclass[12pt]{amsart}
\usepackage{amssymb}

\newcommand{\A}{\mathcal{A}}
\newcommand{\B}{\mathcal{B}}
\newcommand{\C}{\mathbb{C}}
\newcommand{\ca}{\text{\textcircled{\smaller[5]$\aleph_0$\larger[5]}}}
\newcommand{\eps}{\varepsilon}
\newcommand{\F}{\mathbb{F}}
\newcommand{\h}{\mathfrak{H}}
\newcommand{\I}{\mathcal{I}}
\newcommand{\kom}{\mathcal{K}}
\newcommand{\M}{\mathcal{M}}
\newcommand{\Mn}{\mathbb{M}_n}
\newcommand{\Mr}{\mathbb{M}_r}
\newcommand{\n}{\text{\textcircled{$n$}}}
\newcommand{\N}{\mathcal{N}}
\newcommand{\p}{\mathcal{P}}
\newcommand{\pM}{(\p(\M)/\sim)}
\newcommand{\R}{\mathcal{R}}

\newcommand{\U}{\mathcal{U}}

\newcommand{\worb}{\overline{\mathcal{U}(x)}^{\sigma-w}}
\newcommand{\z}{\mathcal{Z}}
\newcommand{\Z}{\mathbb{Z}}

\newtheorem{theorem}{Theorem}
\newtheorem{proposition}[theorem]{Proposition}
\newtheorem{corollary}[theorem]{Corollary}
\newtheorem{lemma}[theorem]{Lemma}

\theoremstyle{definition}
\newtheorem{definition}[theorem]{Definition}

\newtheorem{question}[theorem]{Question}
\newtheorem{problem}[theorem]{Problem}
\newtheorem*{notation}{Notation}

\theoremstyle{remark}
\newtheorem{remark}[theorem]{Remark}

\numberwithin{equation}{section}
\numberwithin{theorem}{section}

\begin{document}

\title[Divisible operators]{Divisible operators in von Neumann algebras}
\author{David Sherman}
\address{Department of Mathematics\\ University of Virginia\\ P.O. Box 400137\\ Charlottesville, VA 22904}
\email{dsherman@virginia.edu}
\subjclass[2000]{Primary 47C15; Secondary 47A65, 46L10}
\keywords{divisible operator, von Neumann algebra, ultrapower, Voiculescu's theorem, approximate equivalence, unitary orbit}

\begin{abstract}
Relativizing an idea from multiplicity theory, we say that an element $x$ of a von Neumann algebra $\M$ is $n$-divisible if $W^*(x)' \cap \M$ unitally contains a factor of type $\text{I}_n$.  We decide the density of the $n$-divisible operators, for various $n$, $\M$, and operator topologies.  The most sensitive case is $\sigma$-strong density in $\text{II}_1$ factors, which is closely related to the McDuff property.

We make use of Voiculescu's noncommutative Weyl-von Neumann theorem to obtain several descriptions of the norm closure of the $n$-divisible operators in $\B(\ell^2)$.  Here are two consequences: (1) in contrast to the reducible operators, of which they form a subset, the divisible operators are nowhere dense; (2) if an operator is a norm limit of divisible operators, it is actually a norm limit of unitary conjugates of a single divisible operator.

This is related to our ongoing work on unitary orbits by the following theorem, which is new even for $\B(\ell^2)$: if an element of a von Neumann algebra belongs to the norm closure of the $\aleph_0$-divisible operators, then the $\sigma$-weak closure of its unitary orbit is convex.
\end{abstract}

\maketitle

\section{Introduction} \label{S:intro}
Let $\B(\h)$ be the algebra of bounded linear operators on a Hilbert space $\h$, and let $x \in \B(\h)$.  The operator $x \oplus x \in \B(\h \oplus \h)$, which applies $x$ to each summand simultaneously, may be thought of as the ``double" of $x$.  Now the latter algebra is just $\mathbb{M}_2(\B(\h))$, the $2 \times 2$ matrices over $\B(\h)$, and it suggests how to double an operator $x$ in an unrepresented von Neumann algebra $\M$: let $x \oplus x \in \mathbb{M}_2(\M)$ be the matrix with $x$ on the diagonal and zeroes elsewhere.  Similarly one may take larger (even infinite) multiples of $x$.  For any cardinal $n$, we employ Ernest's notation (\cite{E}) and set
\begin{equation} \label{E:nx}
\n x \overset{\text{def}}{=} 1 \otimes x \in \Mn \overline{\otimes} \M.
\end{equation}
Here we write $\Mn$ for the factor of type $\text{I}_n$, even when $n$ is infinite.  Note that a multiple of $x \in \M$ belongs to an algebra which may not be isomorphic to $\M$.

We will want to know when a given $x \in \M$ can be written as $\n y$ for some $y$.  In other words, when are there an algebra $\N$ and an isomorphism
\begin{equation} \label{E:div?}
\pi: \M \overset{\sim}{\to} \Mn \overline{\otimes} \N, \qquad \pi(x) \in 1 \otimes \N?
\end{equation}
Clearly this would imply that the relative commutant $W^*(x)' \cap \M$ unitally contains $\Mn$.  And the converse is also valid, since the existence of such an $\Mn$ guarantees an ``internal" isomorphism
\begin{equation} \label{E:inttensor}
\pi: \M \overset{\sim}{\to} \Mn \overline{\otimes} (\Mn' \cap \M), \qquad \pi(x) \in 1 \otimes (\Mn' \cap \M)
\end{equation}
(\cite[Lemma 6.6.3 and Example 11.2.2]{KR}).  We therefore make the following

\begin{definition} \label{T:def}
Let $n$ be a cardinal greater than 1, and let $\M$ be a von Neumann algebra.  For $x \in \M$, we say that $x$ is \textbf{$n$-divisible} if the relative commutant $W^*(x)' \cap \M$ unitally contains $\Mn$.  Similarly, for a $C^*$-algebra $\A$, we say that a *-homomorphism $\rho: \A \to \M$ is \textbf{$n$-divisible} if the relative commutant $\rho(\A)' \cap \M$ unitally contains $\Mn$.  We also say that such an $x$ or $\rho$ is \textbf{divisible} if it is $n$-divisible for some $n$.
\end{definition}

One visualizes such an $x$ as $\left( \begin{smallmatrix} y & 0 & \hdots \\ 0 & y & \hdots \\ \vdots & \vdots & \ddots \end{smallmatrix} \right)$; these are the operators which can be ``divided by $n$."  But there is no hope of defining an operator quotient: if $\pi$ solves \eqref{E:div?}, so does $(\text{id} \otimes \alpha) \circ \pi$ for any $\alpha \in \text{Aut}(\N)$.  Division, unlike multiplication, necessarily involves isomorphism classes.  We give further discussion of this in an Appendix.


We will repeatedly use the following consequence of the basic structure theory of von Neumann algebras.  The unfamiliar reader may see \cite[Lemma 6.5.6 and Section 6.6]{KR}.

\begin{lemma} \label{T:23div}
For $2 \leq n < \aleph_0$, a von Neumann algebra unitally contains $\Mn$ unless it has a type $\text{I}_k$ summand for some $k$ not divisible by $n$.  A von Neumann algebra is properly infinite if and only if it contains $\B(\ell^2)$ unitally.
\end{lemma}

Divisibility is really just a variant of multiplicity.  Recall that a Hilbert space operator $x$ is said to have \textit{(uniform) multiplicity $n$} (or be \textit{homogeneous of order $n$}) if the commutant $W^*(x)'$ is a type $\text{I}_n$ algebra.  See \cite{Hal1951,K1957,A1976,E} for some characterizations, extensions to representations, and applications, especially in regard to the problem of unitary equivalence.  While not every Hilbert space operator has a multiplicity, the type decomposition of the von Neumann algebra $W^*(x)'$ allows us to write $x$ as a direct sum of operators which do, plus an additional term whose commutant has no type I summand.  This readily generalizes to operators in von Neumann algebras by considering the relative commutant $W^*(x)' \cap \M$.  (One loss is that $W^*(x)$ and $W^*(x)' \cap \M$ may have different types.  Multiplicity theory solves the problem of unitary equivalence for normal operators in $\B(\h)$, but it is insufficient for analogous questions in general von Neumann algebras.  See \cite[Section 8]{S2005} for more discussion of this, or \cite{Bu} for a study of multiplicity as an (incomplete) invariant for abelian subalgebras of von Neumann algebras.)  At least for $n$ finite, then, $n$-divisibility of $x$ is equivalent to saying that $W^*(x)' \cap \M$ lacks a $\text{I}_k$ summand for all $k$ indivisible by $n$.

\smallskip

Let us consider some basic questions about the set of $n$-divisible operators in $\M$.  \textit{Can it be empty?}  Yes, if $\M$ has a nonzero type $\text{I}_k$ summand for some cardinal $k$ which is not a multiple of $n$ (by Lemma \ref{T:23div}).  The converse is also true, since otherwise the identity is $n$-divisible.  \textit{Can it be all of $\M$?}  It cannot if $\M$ has separable predual, because then any maximal abelian *-subalgebra is generated by a single operator (\cite{vN}), and such an operator is not $n$-divisible.  But there are large von Neumann algebras in which all elements are $n$-divisible; examples arise in Corollary \ref{T:all}.  \textit{When is it dense in $\M$ with respect to the various topologies?}  This last question is the focus of Sections \ref{S:op} and \ref{S:norm}.  See Table 1.

\begin{table}[b]
\begin{center}
\begin{tabular}{|c|c|c|c|}
\hline
Type & norm & $\sigma$-strong or $\sigma$-strong* & $\sigma$-weak \\ \hline
$\text{I}_k, \: k < \aleph_0$ & N & N & sometimes/N \\ \hline
$\text{I}_\infty$ & N & Y & Y \\ \hline
$\text{II}_1$ & ?/N & sometimes/N & Y/N \\ \hline
$\text{II}_\infty$ & ?/N & Y & Y \\ \hline
III & ? & Y & Y \\ \hline
\end{tabular} \caption{Density of the $n$-divisible operators for different topologies and types of von Neumann algebras, $2 \leq n \leq \aleph_0$.  Where two answers are present, the first refers to $n$ finite and the second to $n = \aleph_0$.}
\end{center}
\end{table}

Recall that a Hilbert space operator is said to be \textit{reducible} if it has a reducing subspace, i.e., if it can be written as $y \oplus z$.  This amounts to requiring that $W^*(x)'$ contain a nontrivial projection, or by the double commutant theorem, that $W^*(x) \neq \B(\h)$.  The eighth of Halmos's ten famous operator theory problems from 1970 (\cite{Hal1970}) asked whether the reducible operators in $\B(\ell^2)$ are norm dense.  The affirmative answer is a consequence of Voiculescu's noncommutative Weyl-von Neumann theorem (\cite{V}, see Theorem \ref{T:voic} below).

Now divisible operators are apparently reducible -- we have slightly upgraded our request of $W^*(x)'$, from nontriviality to the containment of matrix units.  Using Voiculescu's theorem, we obtain several descriptions of the closures of the $n$-divisible operators in $\B(\ell^2)$ (Theorems \ref{T:appvoic} and \ref{T:appvoic2}).  An interesting consequence is that a Hilbert space operator which is a norm limit of divisible operators is actually a norm limit of unitary conjugates of a single divisible operator (Theorem \ref{T:ue}).  Finally we show that the set of all divisible operators in $\B(\ell^2)$ is nowhere dense in the norm topology (Theorem \ref{T:nowhere}).  We derive this from more general statements, but its essence is simple: in any open ball in $\B(\ell^2)$, there is an element $y$ such that $C^*(y)$ contains a rank one projection; this projection is not a norm limit of divisible operators, so $y$ cannot be either.

Unfortunately we have made little headway on similar problems in other von Neumann algebras, where we lack an adequate substitute for Voiculescu's theorem.  (Giordano and Ng recently proved that a certain form of Voiculescu's theorem is true only in injective algebras (\cite{G}).  This interesting result has not yet appeared in final form, and it seems not to have any direct implications for the questions of the present paper.)  But if we work instead with the $\sigma$-strong*, $\sigma$-strong, or $\sigma$-weak topologies, our answers are rather complete and mostly determined by the type of the algebra.  The most delicate case is $\sigma$-strong density in a $\text{II}_1$ algebra, which is closely related to the McDuff property (Theorem \ref{T:mcduff}).  This leads to a fairly simple (but poorly understood) numerical invariant which measures the ``McDuffness" of a singly-generated $\text{II}_1$ factor (Remark \ref{T:mcdrem}).

\smallskip

Seemingly unrelated is

\begin{question} \label{T:conj}
When is the $\sigma$-weakly closed unitary orbit of an element in a von Neumann algebra convex?
\end{question}

For $\B(\h)$, the answer is frequently no, even for self-adjoint operators (\cite{AS}).  But nonatomic factors exhibit different behavior: the answer is yes for all self-adjoint operators, and we know no operator for which the answer is no.  One may view an affirmative answer as a noncommutative Lyapunov-type theorem (\cite{AA}), so it is not surprising that noncommutative atomic measures are recalcitrant.

We show in Section \ref{S:convex} that Question \ref{T:conj} has an affirmative answer whenever the element belongs to the norm closure of the $\aleph_0$-divisible operators (Theorem \ref{T:conv}).  Since the $\aleph_0$-divisible operators are not norm dense in a von Neumann algebra with a semifinite summand (Corollary \ref{T:semif}), this does not give a full answer to Question \ref{T:conj}.  However, at present we do not know if they are dense in a type III algebra; if so, Question \ref{T:conj} has an affirmative answer for all operators in type III algebras.  In fact this was our motivation for studying divisible operators in the first place.

\smallskip

Let us review some of our assumptions and notations.  All $C^*$-algebras, *-homomorphisms, and inclusions are assumed to be unital.  Generic von Neumann algebras are denoted by $\M$ and $\N$, and like the $C^*$-algebras of this paper, they are not assumed to be represented on a Hilbert space.  For this reason we use only intrinsic topologies: the norm, $\sigma$-strong*, $\sigma$-strong, and $\sigma$-weak topologies are symbolized by $\|\|$, $\sigma-s^*$, $\sigma-s$, and $\sigma-w$ respectively.  (We frequently use that the $\sigma$-strong and $\sigma$-strong* topologies agree on finite algebras (\cite[Exercise V.2.5]{T}).)  We write $\z(\M)$ for the center and $\U(\M)$ for the unitary group.  As already mentioned, $\Mn$ stands for a factor of type $\text{I}_n$, but we prefer to write $\B(\ell^2)$ for $\mathbb{M}_{\aleph_0}$.  We let $\F_n$ be the free group on $n$ generators, and $L(G)$ be the von Neumann algebra generated by the left regular representation of the group $G$.  The hyperfinite $\text{II}_1$ factor is denoted by $\R$.

We write the normalized trace on a finite type I factor as ``Tr".  For any faithful normal tracial state $\tau$ on a finite algebra $\M$, the Hilbert space $L^2(\M, \tau)$ is obtained by endowing $\M$ with the inner product $\langle x \mid y \rangle = \tau(y^*x)$ and completing in the induced norm $\|x\|_2 = \sqrt{\tau(x^*x)}$.  It is easy to check that $L^2(\M, \tau)$ is a contractive $\M-\M$ bimodule: $\|xyz\|_2 \leq \|x\| \|y\|_2 \|z\|$.  This leads to the useful fact that on bounded subsets of $\M$, the $L^2$ norm determines the $\sigma$-strong topology (\cite[Proposition III.2.2.17]{Bl}).

For $x \in \M$, the left (resp. right) support projection $s_\ell(x)$ (resp. $s_r(x)$) is defined to be the support of $xx^*$ (resp. $x^*x$).  The unitary orbit is $\U(x) = \{uxu^* \mid u \in \U(\M)\}$.  If $x$ is normal and $E \subseteq \C$ is Borel, the corresponding spectral projection is $\chi_E(x)$.

An isomorphism between operators in algebras $(x \in \M) \cong (y \in \N)$ means a *-isomorphism of algebras taking one operator to the other.  The algebras may be omitted when they are understood or irrelevant.  For $\M \cong \N \cong \B(\h)$, this is unitary equivalence, but in general it is a weaker equivalence relation.

The following lemma collects some basic observations about divisibility.

\begin{lemma} \label{T:facts}
Let $n$, $p$ and $r$ be cardinals satisfying $np = r$, and let $x,y,z$ be elements of von Neumann algebras, with $z$ a central projection in the same algebra as $x$.
\begin{enumerate}
\item $x \cong y \Rightarrow \n x \cong \n y$.
\item $\n (\text{\textcircled{$p$}} x) \cong \text{\textcircled{$r$}} x$.  (So $r$-divisibility implies $n$-divisibility.)
\item If $x$ is $n$-divisible, so is every element of $W^*(x)$.
\item In any of the topologies under discussion, if $x$ is a limit of $n$-divisible (or divisible) operators, so is $zx$.
\item If $x$ is a norm limit of $n$-divisible (or divisible) operators, so is every element of $C^*(x)$.
\item If $x$ is a $\sigma$-strong* limit of $n$-divisible (or divisible) operators, then it is a $\sigma$-strong* limit of a uniformly bounded net of $n$-divisible (or divisible) operators.
\item If $x$ is a $\sigma$-strong* limit of $n$-divisible (or divisible) operators, so is every element of $W^*(x)$.
\end{enumerate}
\end{lemma}

\begin{proof}
We prove only the last two statements.

Part (6) goes by the technique of the Kaplansky density theorem, as follows.  Suppose $\{x_\alpha\}$ converges $\sigma$-strong* to $x \in \M$, with $x_\alpha$ $n_\alpha$-divisible and $\{e^\alpha_{ij}\}_{1 \leq i,j \leq n_\alpha}$ matrix units in the relative commutant.  Then $\{(\begin{smallmatrix} 0 & x_\alpha \\ x_\alpha^* & 0 \end{smallmatrix} )\}$ converges $\sigma$-strongly to $(\begin{smallmatrix} 0 & x \\ x^* & 0 \end{smallmatrix} )$.  Consider the continuous truncation function
$$f: \mathbb{R} \to [-\|x\|, \|x\|]; \qquad t \mapsto \begin{cases} -\|x\|, & t \leq -\|x\|; \\ t, & -\|x\| \leq t \leq \|x\|; \\ \|x\|, & \|x\| \leq t. \end{cases}$$
Being continuous and bounded, this function is $\sigma$-strongly continuous from $\M_{sa}$ to the self-adjoint part of the closed ball of radius $\|x\|$, which it fixes pointwise.  (This is already contained in Kaplansky's original paper \cite{Ka1951}; note that the strong and $\sigma$-strong topologies agree for a suitable representation of $\M$.)  Thus
$$f((\begin{smallmatrix} 0 & x_\alpha \\ x_\alpha^* & 0 \end{smallmatrix})) \overset{\sigma-s}{\to} f((\begin{smallmatrix} 0 & x \\ x^* & 0 \end{smallmatrix})) = (\begin{smallmatrix} 0 & x \\ x^* & 0 \end{smallmatrix}).$$
This implies that the 1,2-entries of $f((\begin{smallmatrix} 0 & x_\alpha \\ x_\alpha^* & 0 \end{smallmatrix}))$ converge $\sigma$-strong* to $x$.  By construction these entries are bounded in norm by $\|x\|$.  To finish the proof it suffices to note that these entries are $n_\alpha$-divisible, and that follows from the calculation
$$\{f((\begin{smallmatrix} 0 & x_\alpha \\ x_\alpha^* & 0 \end{smallmatrix}))\}' \supseteq \{(\begin{smallmatrix} 0 & x_\alpha \\ x_\alpha^* & 0 \end{smallmatrix})\}' \supseteq \{I_2 \otimes e_{ij}^\alpha\}.$$

For (7), suppose that $\{x_\alpha\}$ is a bounded net of divisible operators converging $\sigma$-strong* to $x$.  Since multiplication is jointly $\sigma$-strong* continuous on bounded sets, for any noncommutative polynomial $p(x,x^*)$ we have that $\{p(x_\alpha, x_\alpha^*)\}$ is a net of divisible operators converging $\sigma$-strong* to $p(x,x^*)$.  But these last are $\sigma$-strong* dense in $W^*(x)$.
\end{proof}

One would not expect that statement (5) can be strengthened to include all elements of $W^*(x)$, and in fact it cannot.  The relevant example is postponed to Remark \ref{T:now}.

\smallskip

We end this Introduction by pointing out some relations to the literature.


In the proofs of Theorem \ref{T:typeI}(1,2), Proposition \ref{T:ultra}, and Theorem \ref{T:ncomp} below we use the fact (valid in the circumstances of these theorems, not in general) that an element is a limit of $n$-divisible operators if and only if it commutes arbitrarily well with some system of $n \times n$ matrix units.  This is mentioned here for two reasons.  First, it suggests central sequences, which are a key tool in Section \ref{S:op}.  And second, it has a variant which has proved useful in the $K$-theory of $C^*$-algebras: a $C^*$-algebra is \textit{approximately divisible} if any finite set of elements commutes arbitrarily well (in norm) with the unit ball of some finite-dimensional $C^*$-subalgebra having no abelian summand (\cite[Definition 1.2]{BKR}).  Our Theorem \ref{T:mcduff} is in the same vein as \cite[Example 4.8]{BKR}.




Halpern (\cite{Halp1977}) considered elements $x \in \M$ for which $W^*(x)' \cap \M$ contains two complementary projections which are equivalent in $\M$, not necessarily equivalent in $W^*(x)' \cap \M$.  This weaker condition is a form of 2-diagonalizability, not 2-divisibility.  Actually, Halpern's condition is satisfied by every self-adjoint operator in a von Neumann algebra lacking finite type $\text{I}_{2k+1}$ summands (\cite[Remark p.134]{Halp1977}).  His results were later extended by Kadison (\cite{K1984}) and Kaftal (\cite{Kaf1991}).  The $C^*$-version is subtler (\cite{GrP}).

The reader may wonder about a connection to the perturbation theory of operator algebras.  The typical setup, as introduced in \cite{KK}, has the unit balls of two operator algebras uniformly close (i.e., small Hausdorff distance); one may often deduce structural similarities between the algebras.  This is too strong a hypothesis for the present paper, but the weaker notion of convergence in the Effros-Mar\'{e}chal topology (\cite{Ef}) has some relevance.  We only wish to mention that in the language of Haagerup-Winsl\o w (\cite[Definition 2.1]{HW1}), $x \in \M$ is a $\sigma$-strong* limit of $n$-divisible operators if and only if $W^*(x)$ is contained in the ``$\liminf$" of a net of subalgebras of $\M$ whose relative commutants contain $\Mn$.

Density questions for divisible operators in \textit{matrices} were carefully studied in a long paper by von Neumann (\cite{vN1942}).  Since this work seems to have gotten little attention for many years and is now being rediscovered by the free probability community, we state a version of its main result: There is $\varepsilon > 0$ such that for every $r \in \mathbb{N}$ there is a contraction $x_r \in \Mr$ which is at least $\varepsilon$ away from the set of divisible operators in the $L^2(\Mr, \text{Tr})$ norm.  In an appropriate context this says that divisible operators do not get arbitrarily dense in the strong topology as the size of the matrix algebra grows.  Other questions raised in \cite{vN1942} still remain open.

\smallskip

\textbf{Acknowledgments.} We are grateful to Chuck Akemann for uncountably many useful comments.  Thanks are also due to Ken Dykema for pointing out some relevant references.

\section{Closures in operator topologies} \label{S:op}

In this section we determine whether the $n$-divisible operators are dense in the operator topologies ($\sigma$-weak, $\sigma$-strong, $\sigma$-strong*), for various algebras and $n$.

\begin{proposition} \label{T:spinf}
In a properly infinite algebra $\M$, the $\aleph_0$-divisible operators are $\sigma$-strong* dense.
\end{proposition}

\begin{proof}
Let $x \in \M$.  Informally, we produce an $\aleph_0$-divisible approximant by taking a big corner of $x$ and copying it down the diagonal $\aleph_0$ times.

Choose a sequence of projections $\{p_n\}$ increasing to 1, with $p_n \sim (1-p_n)$.  Each $p_n$ is a minimal projection in some copy of $\B(\ell^2)$ unitally included in $\M$, so as in \eqref{E:inttensor} there is an isomorphism $\theta_n$ from $\B(\ell^2) \overline{\otimes} p_n \M p_n$ to $\M$ such that
$$\theta_n(e_{11} \otimes y) = y, \qquad \forall y \in p_n\M p_n \subset \M.$$
Note that
$$\theta_n(e_{11}^\perp \otimes p_n) = p_n^\perp \overset{\sigma-s^*}{\to} 0.$$

Now $\theta_n(1 \otimes p_n x p_n)$ is $\aleph_0$-divisible, since $1 \otimes p_n x p_n$ is.  We compute
\begin{align*}
x - \theta_n (1 \otimes p_n x p_n) &= x - \theta_n(e_{11} \otimes p_n x p_n) - \theta_n (e_{11}^\perp \otimes p_n x p_n)) \\ &= x - p_n x p_n - \theta_n (e_{11}^\perp \otimes p_n x p_n) \\ &= p_n^\perp x  + x p_n^\perp - p_n^\perp x p_n^\perp  - p_n^\perp \theta_n (1 \otimes p_n x p_n) p_n^\perp,
\end{align*}
and each term in the last expression converges $\sigma$-strong* to zero.
\end{proof}

\begin{lemma} \label{T:normal}
All normal operators in a type II (resp. III) algebra belong to the norm closure of the $n$-divisible operators, for any $n < \aleph_0$ (resp. $n \leq \aleph_0$).
\end{lemma}

\begin{proof}
We say that $x \in \M$ is a \textit{simple} operator if $x = \sum_{j=1}^n \lambda_j p_j$, where $\{\lambda_j\}$ are distinct scalars and $\{p_j\}$ are projections adding to 1.  In this case it is easy to check that $W^*(x)' \cap \M = \oplus p_j \M p_j$.  It then follows from Lemma \ref{T:23div} that a simple operator in a type II (resp. III) algebra is $n$-divisible for any $n < \aleph_0$ (resp. $n \leq \aleph_0$).

The spectral theorem guarantees that any normal operator is well-approximated in norm by simple operators, finishing the proof.
\end{proof}

\begin{proposition}
For $n < \aleph_0$, the $n$-divisible operators are $\sigma$-weakly dense in any $\text{II}_1$ von Neumann algebra.
\end{proposition}

\begin{proof}
From Lemma \ref{T:normal}, we know that any normal operator is norm-approximated by $n$-divisible operators.  The conclusion follows from the fact that unitaries are $\sigma$-weakly dense in the unit ball of any nonatomic von Neumann algebra (\cite[Theorem 1]{Dy}).
\end{proof}

Halmos (\cite[Proposition 1]{Hal1968}) gave an elementary proof that the reducible operators are not dense in a finite type I factor, so the smaller set of $n$-divisible operators is not dense either.  Of course this also follows from the difficult von Neumann result stated at the end of the Introduction.  The next theorem (plus Proposition \ref{T:spinf}) handles all density questions for type I algebras.  

\begin{theorem} \label{T:typeI}
Let $2 \leq n < \aleph_0$.
\begin{enumerate}
\item The $n$-divisible operators are not $\sigma$-strongly dense in a finite type I algebra.
\item The $n$-divisible operators are not norm dense in a type I algebra.
\item The $n$-divisible operators are $\sigma$-weakly dense in a finite type I algebra if and only if it is nonatomic and only has type $\text{I}_k$ summands for $k$ in some subset of $\{2n, 3n, 4n, \dots\}$.
\end{enumerate}
\end{theorem}

\begin{proof}
(1) Given a finite type I algebra $\M$, let $\tau$ be a normal tracial state with (central) support $z$.  Recall that the $L^2$-norm $\|x\|_2 = \sqrt{\tau(x^* x)}$ determines the $\sigma$-strong topology on bounded subsets of $z\M$.  In the next paragraph we work entirely inside $z\M$ and show that the $n$-divisibles are not $\sigma$-strongly dense, which is sufficient by Lemma \ref{T:facts}(4).

Let $p$ be a nonzero abelian projection, and suppose $\{x_\alpha\}$ were $n$-divisibles converging $\sigma$-strongly to $p$.  By Lemma \ref{T:facts}(6) we may assume that $\{x_\alpha\}$ are uniformly bounded.  Let $\{e^\alpha_{ij}\}_{1 \leq i,j \leq n}$ be matrix units commuting with $x_\alpha$.  Then
$$\|[p, e^\alpha_{ij}]\|_2 = \|[p - x_\alpha, e^\alpha_{ij}]\|_2 \leq 2\|p - x_\alpha\|_2 \to 0$$
for any fixed $i$ and $j$.  Using the fact that $p$ is abelian, we compute for any $i \neq j$,
\begin{align*}
pe^\alpha_{ii}p &= p e^\alpha_{ij} e^\alpha_{jj} e^\alpha_{ji} p \\ &\sim (p e^\alpha_{ij} p)(p e^\alpha_{jj} p)(p e^\alpha_{ji} p) \\ &= (p e^\alpha_{ij} p)(p e^\alpha_{ji} p)(p e^\alpha_{jj} p) \\ &\sim p e^\alpha_{ij} e^\alpha_{ji} e^\alpha_{jj} p \\ &= 0.
\end{align*}
Here ``$\sim$" represents an $L^2$ approximation which gets better as $\alpha$ increases; we conclude $pe^\alpha_{ii}p \to 0$ $\sigma$-strongly for any $i$.  But then
$$p = p \left(\sum_{i=1}^n e^\alpha_{ii} \right) p = \sum_{i=1}^n (pe^\alpha_{ii}p) \overset{\sigma-s}{\to} 0,$$
a contradiction.

(2) The argument in the second paragraph of the proof of (1) also establishes (2), except that one uses the uniform norm instead of the $L^2$/$\sigma$-strong topology.

(3) For the necessity we again appeal to Lemma \ref{T:facts}(4), as follows.  An algebra which fails to be nonatomic has a finite matrix algebra as a direct summand, and as mentioned before the theorem, the $n$-divisible operators are not dense in this summand.  In a summand of type $\text{I}_k$, $k$ indivisible by $n$, there are no $n$-divisible operators at all.  In a type $\text{I}_n$ summand, the $n$-divisible operators coincide with the center, which is already $\sigma$-weakly closed.  It remains to establish sufficiency, i.e. $\sigma$-weak density of the $n$-divisibles in $\mathbb{M}_{mn} \otimes L^\infty(X, \mu)$ for any $m \geq 2$ and nonatomic $(X, \mu)$.

We explain some reductions.  It will be enough to show that any simple function $\sum (x_j \otimes \chi_{E_j})$ is a $\sigma$-weak limit of $n$-divisibles, as the simple functions are $\sigma$-weakly (even norm) dense.  It will also be enough to show that any $x \otimes 1$ is a $\sigma$-weak limit, since this argument can then be applied to the sub-measure spaces $(E_j, \mu\mid_{E_j})$ occurring in a simple function.  Finally note that $L^\infty(X, \mu)$ can be written as the direct sum of $L^\infty$ algebras based on nonatomic probability spaces.  (For details, see \cite[Proof of (iii) $\Rightarrow$ (i) in Theorem III.1.18]{T}.  Any finite measure can be rescaled to be a probability measure without changing the associated $L^\infty$ algebra.)  It then suffices to work in $\N = \mathbb{M}_{mn} \otimes L^\infty(X, \mu)$, where $(X, \mu)$ is a nonatomic probability space.



Next we claim that the linear span of the $n$-divisibles in $\mathbb{M}_{mn}$ is all of $\mathbb{M}_{mn}$.  There are many ways to do this; here is a relatively short one.  Let $p$ be a projection of rank $n$, and consider the linear span $L$ of the unitary orbit $\U(p)$.  According to a result of Marcoux and Murphy (\cite[Theorem 3.3]{MM}), such a unitarily-invariant linear space is either contained in the center or contains all commutators.  Since $p$ is not central, $L$ contains all commutators.  Now write any $y \in \mathbb{M}_{mn}$ as $(y - \text{Tr}(y)1) + \text{Tr}(y)1$.  The first summand has trace zero and so is a commutator (\cite{Sho}).  The second summand is a multiple of the identity, which is a sum of $m$ orthogonal conjugates of $p$.  We conclude that both, and also $y$, belong to $L$.  This establishes that $L= \mathbb{M}_{mn}$.  But $L$ is contained in the linear span of the $n$-divisibles, which is therefore also all of $\mathbb{M}_{mn}$.

We return to the goal of finding $n$-divisibles in $\N$ converging $\sigma$-weakly to $x \otimes 1$.  Since $n$-divisibles are closed under scalar multiplication, by the preceding paragraph we can write $x$ as the average of some finite set of $n$-divisibles: $x = \frac{1}{J}\sum_{j=0}^{J-1} x_j$.  The idea of our remaining step is to ``spread the $x_j$ evenly" over $(X, \mu)$ so as to converge weakly to their average.  For this we need to divide $X$ up into finer and finer pieces.  Using the symbol ``$\sqcup$" for disjoint union, recursively find subsets of $X$ indexed by finite strings in $\{0, 1, \dots, J-1\}$, as follows:
$$E_\varnothing = X;$$
$$E_{j_1 j_2 \dots j_{\ell-1}} = E_{j_1 j_2 \dots j_{\ell-1} 0} \sqcup E_{j_1 j_2 \dots j_{\ell-1} 1} \sqcup \dots E_{j_1 j_2 \dots j_{\ell-1} (J-1)},$$
$$\mu(E_{j_1 j_2 \dots j_{\ell-1} j}) = J^{-\ell} \quad (1 \leq j \leq J).$$
Define
$$y_\ell = \sum_{j=0}^{J-1} x_j \otimes \chi_{H^\ell_j},$$
where $H^\ell_j$ is the disjoint union of all sets $E$ labeled by strings of length $\ell$ ending in $j$.  (To illustrate the idea, assume that $(X, \mu)$ is $([0,1], m)$.  Then one may take $E_{j_1 j_2 \dots j_\ell}$ to be the set of numbers whose expansion in base $J$ begins with the string $j_1 j_2 \dots j_\ell$, so that $H^\ell_j$ is the set of numbers whose $\ell$th digit is $j$.  Numbers in $[0,1]$ with nonunique expansion constitute a null set and may be distributed arbitrarily.)  Note that $y_\ell$ is $n$-divisible.  We assert that $\{y_\ell\}$ converges $\sigma$-weakly to $x \otimes 1$.

To prove this assertion, consider the trace $\tau = \text{Tr} \otimes (\int \cdot \: d\mu)$.  In $L^2(\N, \tau)$ one can verify the following computations, assuming $\ell \neq \ell'$:
$$\langle y_\ell , y_\ell' \rangle = \langle y_\ell , x \otimes 1 \rangle = \langle x \otimes 1 , x \otimes 1 \rangle = \frac{1}{J^2} \sum_{i,j} \text{Tr}(x_i^* x_j);$$
$$\langle y_\ell , y_\ell \rangle = \frac{1}{J} \sum_j \text{Tr}(x_j^* x_j).$$
It follows from these that $\{y_\ell - (x \otimes 1)\}$ is an orthogonal sequence of constant norm and thus converges weakly to 0 in the Hilbert space $L^2(\N, \tau)$.  In particular, for $h \in \N \subset L^2(\N, \tau)$,
\begin{equation} \label{E:L1}
\langle y_\ell - (x \otimes 1), h \rangle \to 0 \quad \text{and so} \quad \langle y_\ell, h \rangle \to \langle x \otimes 1, h \rangle.
\end{equation}
Since functionals of the form $\{\langle \cdot, h \rangle \mid h \in \N\}$ are norm dense in $\N_*$ (\cite[Theorem V.2.18]{T}), \eqref{E:L1} establishes the $\sigma$-weak convergence $y_\ell \to x \otimes 1$, finishing the proof.
\end{proof}




In the rest of this section $\M$ is a $\text{II}_1$ factor with separable predual and trace $\tau$, $2 \leq n < \aleph_0$, and $\omega$ is a free ultrafilter on $\mathbb{N}$.

We recall the construction of the \textit{tracial ultrapower} $\M^\omega$, essentially due to Wright (\cite[Theorems 2.6 and 4.1]{W}; see \cite{S2008} for comments on the historical record).  Let $\I_\omega \subset \ell^\infty(\M)$ be the two-sided ideal of sequences $(x_k)$ with $x_k$ converging $\sigma$-strongly to 0 as $k \to \omega$.  We define $\M^\omega$ to be the quotient $\ell^\infty(\M)/\I_\omega$, shown to be a $\text{II}_1$ factor by Sakai (\cite[Section II.7]{Sa}).  Assuming the continuum hypothesis, $\M^\omega$ and the inclusion $\pi: \M \hookrightarrow \M^\omega$ (as cosets of constant sequences) do not depend on the choice of $\omega$ (\cite[Theorem 3.2]{GH}).





\begin{proposition} \label{T:ultra}
For $x \in \M$, the following conditions are equivalent:
\begin{enumerate}
\item $x$ is a $\sigma$-strong limit of $n$-divisible operators;
\item $\pi(x) \in \M^\omega$ is $n$-divisible.
\end{enumerate}
\end{proposition}

\begin{proof} ${}$
(1) $\Rightarrow$ (2): Appealing to Lemma \ref{T:facts}(6) and the $L^2(\M, \tau)$ topology, there is a sequence $\{x_k\}$ of $n$-divisible operators converging $\sigma$-strongly to $x$.  Then each $x_k$ commutes with a system of matrix units $\{e^k_{ij}\}_{1 \leq i,j \leq n}$.  It is immediate that $\{(e^k_{ij})_k\}_{1 \leq i,j \leq n}$ are matrix units in $\M^\omega$.  A straightforward computation shows that they commute with $\pi(x)$:
\begin{equation} \label{E:matunits}
\|[x, e^k_{ij}]\|_2 = \|[x - x_k, e^k_{ij}]\|_2 \leq 2 \|e^k_{ij}\| \|x - x_k\|_2 \overset{k \to \omega}{\longrightarrow} 0.
\end{equation}

(2) $\Rightarrow$ (1): Suppose that $\pi(x)$ commutes with matrix units in $\M^\omega$, and let the $ij$-unit have representing sequence $(f^k_{ij})_k$.  For each $k$, the set $\{f^k_{ij}\}_{1 \leq i,j \leq n}$ need not consist of matrix units in $\M$, but one may apply an argument of McDuff (\cite[Lemma 8]{McD}) to find true matrix units $\{e^k_{ij}\}_{1 \leq i,j \leq n}$ such that $(f^k_{ij})_k = (e^k_{ij})_k$ for each $1 \leq i,j \leq n$.

Now set
\begin{equation} \label{E:csapprox}
x_k = \sum_{m=1}^n e^k_{m1} x e^k_{1m}.
\end{equation}
It is easy to check that $x_k$ converges to $x$ $\sigma$-strongly as $k \to \omega$:
$$\|x - x_k\|_2 = \|\sum x e^k_{m1} e^k_{1m} - \sum e^k_{m1} x e^k_{1m}\|_2 = \|\sum [x, e^k_{m1}] e^k_{1m}\|_2 \to 0.$$
(Compare \cite[Proof of Lemma 1]{McD}.)  By construction $x_k$ commutes with $\{e^k_{ij}\}_{1 \leq i,j \leq n}$ and so is $n$-divisible.
\end{proof}

The next result produces large von Neumann algebras in which every element is $n$-divisible, as promised in the Introduction.  (The condition on $\M$ is neither universal nor impossible, as we will see momentarily.)

\begin{corollary} \label{T:all}
The following conditions are equivalent for $\M$:
\begin{enumerate}
\item the $n$-divisible operators are $\sigma$-strongly dense in $\M$;
\item every element of $\pi(\M)$ is $n$-divisible in $\M^\omega$;
\item every element of $\M^\omega$ is $n$-divisible.
\end{enumerate}
\end{corollary}

\begin{proof}
Proposition \ref{T:ultra} gives the equivalence of (1) and (2), and (3) is clearly stronger than (2).  So let us assume (1) and show (3).

Let $(x_k)$ represent an element of $\M^\omega$.  For a fixed $k$, $x_k$ is a $\sigma$-strong limit of $n$-divisible operators, so as in \eqref{E:matunits} we may find matrix units $\{e_{ij}^k\}$ with
$$\|[x_k, e_{ij}^k]\|_2 < \frac1k, \qquad \forall 1 \leq i,j \leq n.$$
Then $\{(e_{ij}^k)\}_{1 \leq i,j \leq n}$ are matrix units in $\M^\omega$ which commute with $(x_k)$.
\end{proof}

We need to review some more terminology.  A \textit{generator} of a von Neumann algebra is $x \in \N$ satisfying $\N = W^*(x)$.  The \textit{generator problem} asks if every von Neumann algebra with separable predual has a generator; the only unresolved cases are certain $\text{II}_1$ factors, including in particular $L(\F_3)$.  A recent survey of the generator problem is \cite{Sh}.

We say that a $\text{II}_1$ factor $\M$ is \textit{McDuff} if $\M \cong \M \overline{\otimes} \R$.  The main result of \cite{McD} can then be formulated as follows.

\begin{theorem} \label{T:mcd} $($\cite[Theorem 3 and Lemma 7]{McD}$)$
$\M$ is McDuff if and only if $\pi(\M)' \cap \M^\omega$ is noncommutative, and in this case $\pi(\M)' \cap \M^\omega$ is type $\text{II}_1$.
\end{theorem}

\begin{theorem} \label{T:mcduff}
Let $2 \leq n < \aleph_0$.  If $\M$ is McDuff, the $n$-divisible operators are $\sigma$-strongly dense.  For singly-generated $\M$, the converse holds as well.
\end{theorem}

\begin{proof}
If $\M$ is McDuff, we know from Theorem \ref{T:mcd} that $\pi(\M)' \cap \M^\omega$ is type $\text{II}_1$.  Then for any $x \in \M$,
$$W^*(\pi(x))' \cap \M^\omega \supseteq \pi(\M)' \cap \M^\omega \supset \Mn,$$
and the conclusion follows from Proposition \ref{T:ultra}.

If $\M = W^*(x)$ and $x$ is a $\sigma$-strong limit of $n$-divisibles, then Proposition \ref{T:ultra} again implies
$$\pi(\M)' \cap \M^\omega = W^*(\pi(x))' \cap \M^\omega \supset \Mn.$$  By Theorem \ref{T:mcd}, the noncommutativity of $\pi(\M)' \cap \M^\omega$ implies that $\M$ is McDuff.
\end{proof}

\begin{remark} \label{T:mart}
Here is another way to see that $n$-divisible operators are $\sigma$-strongly dense in a McDuff factor $\M$.  Let $\M_k = \M \otimes \mathbb{M}_{2^k} \subset \M \overline{\otimes} \R \cong \M$ be a sequence of increasing subfactors with $\sigma$-strongly dense union.  Note that each $\M_k$ has relative commutant $\cong \R$, so they consist entirely of $n$-divisible operators.  Let $E_k$ be the trace-preserving conditional expectation from $\M$ onto $\M_k$.  A simple martingale theorem (first proved in \cite[Corollary 2.1]{U}) shows that $E_k(x) \overset{\sigma-s}{\to} x$, for any $x \in \M$.
\end{remark}

The reader will admit the existence of McDuff factors: tensor any finite factor with $\R$.  It may be less clear that there are $\text{II}_1$ factors in which the $n$-divisibles ($2 \leq n < \aleph_0$) are not $\sigma$-strongly dense, so we now provide a variety of examples.  Note that this is \textit{not} intended as a complete list.  The reader is referred to the sources for explanation of undefined terms.

\begin{corollary} \label{T:sometimes}
For any $n > 1$, the $n$-divisible operators are not $\sigma$-strongly dense in any of the following $\text{II}_1$ factors:
\begin{enumerate}
\item $L(SL(k,\Z))$ ($k \geq 3$ and odd) and $L(PSL(k,\Z))$ ($k \geq 4$ and even);
\item tensor products of two $\text{II}_1$ factors, neither McDuff and one with property $T$;
\item $L^\infty(X, \mu) \rtimes G$, where $(X, \mu)$ is a nonatomic probability space, $G$ is a countable discrete non-inner amenable group, and the action is free, ergodic, and measure-preserving;
\item factors which have $\Gamma$ but are not McDuff;
\item $L(\F_m)$ ($m \geq 2$).
\end{enumerate}
\end{corollary}

\begin{proof} For the first four classes, we simply explain why the factor is singly-generated and not McDuff.  The conclusion then follows from Theorem \ref{T:mcduff}.

(1): They have property $T$ and so are not McDuff by \cite{C1980,CJ}.  They are singly-generated by \cite{GS}.

(2): They are not McDuff by \cite[Corollary 3.7]{M}.  Any tensor product of $\text{II}_1$ factors is singly-generated by \cite{GP}.

(3): They are not McDuff by \cite[Proposition 3.9]{M}.  Any $\text{II}_1$ factor with a Cartan subalgebra is singly-generated by \cite{P1985}.

(4): Factors with $\Gamma$ are singly-generated by \cite{GP}.  The first example of a non-McDuff factor with $\Gamma$ was constructed in \cite[Proposition 22]{DL}.

(5): $L(\F_m)$ is not known to be singly-generated for $m \geq 3$, so we use Proposition \ref{T:ultra} instead.  Let $\F_m$ have generators $\{g_j\}_{j=1}^m$, so that $L(\F_m)$ has generators $\{\lambda_{g_j}\}$.  Murray and von Neumann showed that these factors do not have $\Gamma$ and so are not isomorphic to $\R$ (\cite[Section VI.6.2]{MvN4}).  Their original estimates can be adapted to show that for a unitary $u \in L(\F_m)$,
\begin{equation} \label{E:free}
\|u - \tau(u)1\|_2 \leq 14\max \{\|[x, \lambda_{g_1}]\|_2, \|[x, \lambda_{g_2}]\|_2\}.
\end{equation}
(See \cite[Equation XIV(3.3)]{T3}.)  This implies that any sequence of unitaries which asymptotically commutes with $\lambda_{g_1}$ and $\lambda_{g_2}$ must be equivalent to a sequence of scalars.  Note that $i \text{Log}(\lambda_{g_j})$ is self-adjoint, and set
$$x = i \text{Log} (\lambda_{g_1}) + \text{Log} (\lambda_{g_2}).$$
We then have
$$W^*(\pi(x))' \cap L(\F_m)^\omega = W^*(\pi(\lambda_{g_1}), \pi(\lambda_{g_2}))' \cap L(\F_m)^\omega = \C \nsupseteq \Mn. \qedhere$$
\end{proof}

Let $2 \leq n < \aleph_0$, and consider the following conditions on a $\text{II}_1$ factor $\M$ with separable predual:

\begin{enumerate}
\item[(1)] $\M$ is McDuff;
\item[(2)] for every singly-generated subalgebra $\N \subseteq \M$, $\pi(\N)' \cap \M^\omega$ is type $\text{II}_1$;
\item[($3_n$)] for every singly-generated subalgebra $\N \subseteq \M$, $\pi(\N)' \cap \M^\omega$ unitally contains $\Mn$.
\end{enumerate}

Each of these conditions implies its successor, and the last is equivalent to $\sigma$-strong density of the $n$-divisible operators in $\M$.  It seems natural to call condition (2) ``locally McDuff."

\begin{problem} \label{T:lmcd}
Is either of the implications (2) $\Rightarrow$ (1), ($3_n$) $\Rightarrow$ (2) valid?
\end{problem}

Both of these implications would follow from an affirmative answer to the generator problem.  This means that to \textit{disprove} one of them, one would have to establish the existence of a von Neumann algebra with separable predual which is not singly-generated.  In posing this problem, we are really asking if either implication can be proved directly, without resolution of the generator problem.

For $\N \subset \M$, the algebra $\pi(\N)' \cap \M^\omega$ has received occasional attention in the literature; see \cite[Lemma 2.6]{C1976}, \cite[Theorem 3.5]{M}, \cite[Lemma 3.3.2]{P1995}, and \cite[Theorems 3.5 and 4.7]{FGL}.  It should not be confused with $\pi(\M)' \cap \N^\omega$, which was studied by Bisch (\cite{B}).  For $\N$ a factor, he showed that the latter algebra is noncommutative exactly when $[\N \subset \M] \cong [\N \overline{\otimes} \R \subset \M \overline{\otimes} \R]$; in this case the inclusion is said to be McDuff.




\begin{remark} \label{T:mcdrem}
Let $\M$ be a singly-generated $\text{II}_1$ factor with separable predual.  By Theorem \ref{T:mcduff}, $\M$ is McDuff if and only if the 2-divisible operators are $\sigma$-strongly dense.  There are a variety of ways to quantify this and so obtain a numerical invariant which measures ``how far $\M$ is from being McDuff."

One way would be to find the supremum of distances (in the $\|\|_2$-metric) from $x$ to the 2-divisible operators, where $x$ runs over the unit ball.  Here is a related approach.  As shown in the proof of Proposition \ref{T:ultra}, if $x$ is strongly approximated by $2$-divisible operators, then approximants may be ``built out of $x$" in the sense of \eqref{E:csapprox}.  So we may instead ask for the distance to operators of the form $v^*v x v^*v + vxv^*$, with $v$ a partial isometry satisfying $v^*v + vv^* = 1$.  This gives the invariant
\begin{equation} \label{E:mcdindex}
\sup_{x \in \M_1} \inf_{\substack{v = v v^* v \\ v^* v + vv^* = 1} } \|x - (vxv^* + v^*v x v^* v)\|_2,
\end{equation}
which is zero if and only if $\M$ is McDuff.  At this point the author knows nothing interesting about this quantity when it is nonzero.  One may define similar invariants based on $n$-divisibles for other $n$; the author also does not know how these numbers depend on $n$.
\end{remark}

\section{Closures in the norm topology} \label{S:norm}
The ostensible goal of this section is to describe the norm closure of the $n$-divisible operators in various von Neumann algebras, but our results are rather incomplete for algebras other than $\B(\h)$.  To some extent this deficiency is due to the lack of a generalization of Voiculescu's theorem -- see \cite{H1998,DH,S2005,S2008} for discussion and partial results.  In $\B(\h)$, at least, we arrive at clean descriptions and ultimately show that the $n$-divisibles are nowhere dense.

A first hope might be to imitate the techniques of the previous section.  There we saw that for singly-generated $\text{II}_1$-factors, $\sigma$-strong density of $n$-divisibles ($2 \leq n < \aleph_0$) is equivalent to the existence of noncommuting central sequences.  Central sequences of matrix units give a ``universal formula" for producing $n$-divisible $\sigma$-strong approximants out of any element, as in \eqref{E:csapprox}.  Can a similar construction work in the norm topology?

The most natural setup is this.  For $\M$ a von Neumann algebra, the quotient $(\ell^\infty(\M)/c_0(\M))$ is a $C^*$-algebra.  Let $\sigma: \M \hookrightarrow (\ell^\infty(\M)/c_0(\M))$ be the inclusion as (cosets of) constant sequences.  Then the ``central sequence algebra" is $\sigma(\M)' \cap (\ell^\infty(\M)/c_0(\M))$.  If it were to contain $\Mn$ unitally, one could mimic \eqref{E:csapprox} and use the matrix units to build $n$-divisible norm approximants for any operator.  Unfortunately, this sort of central sequence algebra is always commutative.

\begin{proposition} \label{T:nocent} Under the circumstances of the preceding paragraph,
$$\sigma(\M)' \cap (\ell^\infty(\M)/c_0(\M)) = \ell^\infty(\z(\M))/c_0(\z(\M)).$$
\end{proposition}

\begin{proof}
If $(x_n)$ represents an element of the left-hand side, then
$$\| (\text{ad }x_n) (y) \| = \|[x_n, y]\| \to 0, \qquad \forall y \in \M.$$
By \cite[Theorem 3.4]{El}, a sequence of derivations converges in the point-norm topology only if it converges in the norm topology, so $\| \text{ad }x_n\| \to 0$.  We also have that $\| \text{ad }x_n\| = 2 \text{ dist}(x_n, \z(\M))$ (proved independently in \cite{Ga} and \cite{Z}).  Therefore $(x_n)$ can also be represented by a sequence from $\z(\M)$.
\end{proof}

\begin{remark}
In contrast to Proposition \ref{T:nocent}, Hadwin's asymptotic double commutant theorem (\cite{H1978b}) implies that for any \textit{separable} subset $S \subset \B(\ell^2)$, $\sigma(S)' \cap (\ell^\infty(\B(\ell^2))/c_0(\B(\ell^2)))$ is noncommutative.
\end{remark}

Since infinite-dimensional von Neumann algebras are not norm separable and therefore never singly-generated (or even countably-generated) as $C^*$-algebras, Proposition \ref{T:nocent} cannot be used to preclude the density of $n$-divisibles in the manner of Theorem \ref{T:mcduff}.

\smallskip

The symbol $\kom$ will denote the closed *-ideal generated by the finite projections in any von Neumann algebra under discussion.

\begin{proposition} \label{T:distance}
Let $\M$ be a properly infinite semifinite von Neumann algebra and $k \in \kom$.  Then the distance from $k$ to the $\aleph_0$-divisible operators is $\frac{\|k\|}{2}$.  In particular, if $k$ is a norm limit of $\aleph_0$-divisible operators, then $k=0$.
\end{proposition}

\begin{proof} Let $y$ be $\aleph_0$-divisible, and fix $\eps > 0$.  Then $p = \chi_{[\|y\| - \eps, \|y\|]} (|y|)$ is also $\aleph_0$-divisible and therefore infinite.   In order to mesh cleanly with the paper \cite{K}, represent $\M$ faithfully on a Hilbert space $\h$.  By \cite[Theorem 1.3(d)]{K}, $k$ is not bounded below on $p \h$, so there is a unit vector $\xi \in p \h$ with $\|k \xi\| < \eps$.  Since $y$ is bounded below on $p \h$ by $\|y\| - \eps$,
$$\|(y-k) \xi \| \geq \| y \xi \| - \|k \xi\| > \|y\| - 2 \eps.$$
Now $\eps$ is arbitrary, so $\|y-k\| \geq \|y\|$.  Then
$$\|y-k\| \geq \|y\| \geq \|k\| - \|y-k\| \quad \Rightarrow \quad \|y-k\| \geq \frac{\|k\|}{2}.$$
This shows that the distance from $k$ to the $\aleph_0$-divisibles is $\geq \frac{\|k\|}{2}$.

For the opposite inequality, take any $\eps > 0$.  By definition, $k$ is approximated within $\eps$ by an operator $f$ whose supports are finite; let $q = s_\ell(f) \vee s_r(f)$ (which is finite).  Then
$$\|k - qkq\| \leq \|k - f\| + \|f - qkq\| = \|k - f\| + \|q(k-f)q\| < 2\eps.$$
By the finiteness of $q$, we can find a projection $r \geq q$ such that $r \sim r^\perp$.  Note that
$$\|k - rkr\| = \|k - qkq + r(qkq - k)r\| \leq \|k - qkq\| + \|r(qkq - k)r\| \leq 4 \eps.$$
Write $1 = \sum_{j=1}^\infty r_j$, with $r = r_1$.  Let $\{v_j\}$ be partial isometries with $v_1 = r_1$ and $v_j v_j^* = r_j$, $v_j^* v_j = r_1$ for $j\geq 2$.  Finally, consider the $\aleph_0$-divisible operator
$\sum \frac{v_j k v_j^*}{2}$.  We have
\begin{align*}
\left\|k - \left(\sum \frac{v_j k v_j^*}{2}\right) \right\| &\leq \|k - rkr\| + \left\|rkr - \left(\sum \frac{v_j k v_j^*}{2}\right)\right\| \\ &\leq 4\eps + \left\|\frac{rkr}{2} - \left(\sum_{j=2}^\infty v_j \left(\frac{rkr}{2}\right) v_j^*\right)\right\| \\ &\overset{(*)}{=} 4\eps + \frac{\|rkr\|}{2} \\ &\leq 4\eps + \frac{\|rkr - k\| + \|k\|}{2} \\ &\leq 6 \eps + \frac{\|k\|}{2}.
\end{align*}
(The equality $(*)$ is justified by noting that the summation in the previous expression is an orthogonal sum of operators unitarily conjugate to $\frac{rkr}{2}$.)
\end{proof}

Since there are no $\aleph_0$-divisible operators in a finite algebra, we deduce

\begin{corollary} \label{T:semif}
The $\aleph_0$-divisible operators are not norm dense in any semifinite algebra.
\end{corollary}

At this point all entries of Table 1 have been justified.  We next obtain much more specific information for infinite type I factors, writing $0_\infty$ for the zero operator on $\ell^2$.

\smallskip

We need to recall Voiculescu's theorem and the relevant terminology.  Two operators $x,y$ are said to be \textit{approximately equivalent} when there is a sequence of unitaries $\{u_n\}$ with $u_n x u_n^* \to y$ in norm.  (Sometimes this term implies also that the differences $u_n x u_n^* - y$ are compact, but we do not make this requirement here.)  Similarly, two nondegenerate representations $\rho, \sigma$ of a $C^*$-algebra $\A$ are approximately equivalent when there is a net of unitaries $\{u_\alpha\}$ with $(\text{Ad }u_\alpha) \circ \rho \to \sigma$ in the point-norm topology.  We denote approximate equivalence by $\sim_\text{a.e.}$.  It is clear that an approximate equivalence can be multiplied by an arbitrary cardinal, in the sense of Lemma \ref{T:facts}(1).  Note that when $\A = C^*(x)$,
\begin{equation} \label{E:singlyae}
\rho \sim_\text{a.e.} \sigma \iff \rho(x) \sim_\text{a.e.} \sigma(x).
\end{equation}

\begin{notation}  Let $\rho: \A \to \B(\h)$ be a representation of a separable $C^*$-algebra on a separable Hilbert space.  The set
$$\I_\rho \overset{\text{def}}{=} \rho^{-1}(\rho(\A) \cap \kom)$$
is an ideal of $\A$.  This allows us to write $\rho = \rho^1 \oplus \rho^2$, where $\rho^1$ is the restriction of $\rho$ to the reducing subspace $\rho(\I_\rho) \h$, sometimes called the \textit{essential part} of $\rho$ (\cite[p.341]{A1977}).  Of course $\rho^1$ or $\rho^2$ may be absent from this decomposition.
\end{notation}

\begin{theorem} \label{T:voic} $($\cite[Theorem 1.5]{V}$)$
Let $\rho_j$ ($j=1,2$) be representations of a separable $C^*$-algebra $\A$ on a separable Hilbert space $\h$.  Then $\rho_1 \sim_\text{a.e.} \rho_2$ if and only if
\begin{enumerate}
\item[(i)] $\ker \rho_1 = \ker \rho_2$,
\item[(ii)] $\I_{\rho_1} = \I_{\rho_2}$, and
\item[(iii)] $\rho_1^1$ and $\rho_2^1$ are unitarily equivalent when restricted to this ideal.
\end{enumerate}
\end{theorem}

\begin{theorem} \label{T:ncomp}
Let $\h$ be an infinite-dimensional Hilbert space, $n < \aleph_0$, and $k \in \kom \subset \B(\h)$.  Then $k$ is a norm limit of $n$-divisibles if and only if $k \oplus 0_\infty$ is $n$-divisible.
\end{theorem}

\begin{proof}
Assume that $k$ is a norm limit of $n$-divisibles, and let $k \oplus 0_\infty = \text{Re}(k \oplus 0_\infty) + i \text{Im}(k \oplus 0_\infty)$ be the decomposition into real and imaginary parts.  Each is self-adjoint and compact, so we may list their (finitely or infinitely many) nonzero eigenvalues as follows, including multiplicity.
$$\text{Re}(k \oplus 0_\infty): \qquad \lambda_{-1} \leq \lambda_{-2} \leq \dots < 0 < \dots \leq \lambda_2 \leq \lambda_1$$
$$\text{Im}(k \oplus 0_\infty): \qquad \mu_{-1} \leq \mu_{-2} \leq \dots < 0 < \dots \leq \mu_2 \leq \mu_1$$
We further set
$$p_j = \chi_{\{\lambda_j\}} (\text{Re}(k \oplus 0_\infty)); \qquad q_j = \chi_{\{\mu_j\}} (\text{Im}(k \oplus 0_\infty)).$$

We also have that $k \oplus 0_\infty$ is a norm limit of $n$-divisibles.  (Just add the summand $0_\infty$ onto the $n$-divisible operators converging to $k$.)  Let $\{e^{(m)}_{ij}\}_{i,j = i}^n$ be matrix units commuting with the $m$th operator in the sequence.  It follows as in \eqref{E:matunits} that $\|[e_{ij}^{(m)}, k \oplus 0_\infty]\| \to 0$.  By repeated use of the triangle equality, $\|[e_{ij}^{(m)}, x]\| \to 0$ for every $x \in C^*(k \oplus 0_\infty)$, in particular the $p_j$ and $q_j$ considered above.

For each $m$, $\{p_1 e^{(m)}_{ij} p_1\}$ is an $n^2$-tuple in the unit ball of the finite dimensional space $p_1 \B(\h) p_1$.  Pick a subsequence (still denoted by $m$) which is convergent.  Because $\|[e_{ij}^{(m)}, p_1]\| \overset{m \to \infty}{\to} 0$,
$$(\lim_m p_1 e^{(m)}_{ij} p_1) (\lim_m p_1 e^{(m)}_{kl} p_1) = \delta_{jk} (\lim_m p_1 e^{(m)}_{il} p_1).$$
Therefore $\{\lim p_1 e^{(m)}_{ij} p_1\}_{i,j =1}^n$ is a set of matrix units in $p_1 \B(\ell^2) p_1$.

Now refine the subsequence so that $\{p_{-1} e^{(m)}_{ij} p_{-1}\}$ is a convergent $m$-tuple in $p_{-1} \B(\h) p_{-1}$ (with matrix units as limits).  One may continue refining for $q_1$, then $q_{-1}$, then $p_2, p_{-2}, q_{2}, q_{-2}$, etc.  Extract a diagonal subsequence, still calling the index $m$.

Let $r$ be the supremum of all $p_j$ and $q_j$, so that $r^\perp$ is the infinite-rank projection onto the nullspace $\ker (k \oplus 0_\infty) \cap \ker (k^* \oplus 0_\infty)$.  For each $i$ and $j$, the strong limit of $r e^{(m)}_{ij} r$ exists by the previous paragraph, and these limits form matrix units for $r \B(\h) r$.  They fix all nonzero eigenspaces of $\text{Re}(k \oplus 0_\infty)$ and $\text{Im}(k \oplus 0_\infty)$, so they commute with $r(k \oplus 0_\infty)r$ and $r(k \oplus 0_\infty)^*r$.  It follows that $r(k \oplus 0_\infty)r$ is $n$-divisible.  Now choose any matrix units for $r^\perp \B(\h \oplus \ell^2) r^\perp$, and add them to the corresponding matrix units for $r \B(\h) r$ constructed above.  This produces matrix units for $\B(\h \oplus \ell^2)$ which commute with $k \oplus 0_\infty$, completing the proof of the forward implication.

The opposite implication is trivial when $\h$ has uncountable dimension, as then $k$ and $k \oplus 0_\infty$ are unitarily equivalent.  (Remember that $k$ is compact.)  For separable $\h$ we claim $k \sim_{\text{a.e.}} k \oplus 0_\infty$, so that $k$ is a norm limit of $n$-divisible unitary conjugates of $k \oplus 0_\infty$.  To prove the claim, let $\sigma$ be the representation of $C^*(k)$ on $\ell^2$ with $\sigma(1) = 1$ and $\sigma(k) = 0$.  Apply Voiculescu's theorem to conclude $\text{id} \sim_{\text{a.e.}} \text{id} \oplus \sigma$ as representations of $C^*(k)$.  Then use \eqref{E:singlyae}.
\end{proof}


\begin{remark}
There is no variation of Theorem \ref{T:ncomp} for nonatomic factors which is both useful and true.

Here are examples which show that the compactness of $k$ is indispensible for both its implications.
\begin{itemize}
\item Take $k$ to be a (noncompact) projection of corank 1 in $\B(\ell^2)$, so that $k \oplus 0_\infty$ is $n$-divisible.  If $k$ were a norm limit of $n$-divisibles, then $1-k$ would be, too.  (Just subtract the approximating sequence from $1$.)  But $1-k$ is compact and is shown not to be a norm limit of $n$-divisibles by the theorem.
\item Take $k \in \B(\ell^2(\Z))$ to be the (noncompact) bilateral shift.  It generates a maximal abelian *-subalgebra, so $k \oplus 0_\infty$ is not $n$-divisible.  But the argument in Proposition \ref{T:normal} shows that $k$ is a norm limit of $n$-divisibles, as $k$ can be approximated by simple operators whose spectral projections are all infinite.
\end{itemize}
\end{remark}

\begin{theorem} \label{T:appvoic}
Let $x \in \B(\ell^2)$, $n < \aleph_0$, and $\mbox{\textnormal{id}}$ be the identity representation of $C^*(x)$.  The following conditions are equivalent:
\begin{enumerate}
\item $x$ is a norm limit of $n$-divisible operators;
\item $C^*(x \oplus 0_\infty) \cap \kom$ consists of $n$-divisible operators;
\item $\mbox{\textnormal{id}}^1$ is $n$-divisible;
\item $x$ is approximately equivalent to an $n$-divisible operator.
\end{enumerate}
\end{theorem}

\begin{proof} (4) $\Rightarrow$ (1): The hypothesis implies $x$ is a norm limit of unitary conjugates of a fixed $n$-divisible operator.

(1) $\Rightarrow$ (2): If $x$ is a norm limit of $n$-divisible operators, the same holds for every element of $C^*(x)$.  Now $C^*(x \oplus 0_\infty) \cap \kom = (C^*(x) \cap \kom) \oplus 0_\infty$, and the conclusion follows from Theorem \ref{T:ncomp}.

(2) $\Rightarrow$ (3): We have
$$C^*(x \oplus 0_\infty) \cap \kom = \text{id}^1(\I_{\text{id}}) \oplus \text{id}^2(0) \oplus 0_\infty.$$
Restricted to $\I_{\text{id}}$, $\text{id}^1$ is a direct sum of irreducible representations, the image of each being isomorphic to $\kom$ or a matrix factor.  Condition (2) says that these representations all occur with multiplicities divisible by $n$.  But an irreducible representation of an ideal uniquely induces an irreducible representation of the ambient algebra (\cite[Theorem 1.3.4]{A1976}).  So on all of $C^*(x)$, $\text{id}^1$ is a direct sum of irreducible representations, each with multiplicity divisible by $n$.

(Always $\text{id}^1(x)$ is a direct sum of irreducible operators, each occurring with finite multiplicity.  Hadwin (\cite{H1977}) calls these irreducible operators the \textit{isolated reducing operator-eigenvalues with finite multiplicity}, written $\Sigma_{00}(x)$.  Condition (3) says exactly that every element of $\Sigma_{00}(x)$ has (finite) multiplicity divisible by $n$.)

(3) $\Rightarrow$ (4): Write $\text{id}^1 = \n \rho$ for some representation $\rho$.  By Voiculescu's theorem,
$$\text{id} = \text{id}^1 \oplus \text{id}^2 = (\n \rho) \oplus \text{id}^2 \sim_{\text{a.e.}} (\n \rho) \oplus (\n \text{id}^2) = \n (\rho \oplus \text{id}^2).$$
Plugging in $x$ as in \eqref{E:singlyae}, $x \sim_{\text{a.e.}} \n (\rho(x) \oplus \text{id}^2(x))$.
\end{proof}

\begin{theorem} \label{T:appvoic2}
Let $x \in \B(\ell^2)$ and $\mbox{\textnormal{id}}$ be the identity representation of $C^*(x)$.  The following conditions are equivalent:
\begin{enumerate}
\item $x$ is a norm limit of $\aleph_0$-divisible operators;
\item $C^*(x) \cap \kom = \{0\}$;
\item $\mbox{\textnormal{id}}^1$ is void;
\item $x$ is approximately equivalent to an $\aleph_0$-divisible operator;
\item $x \sim_{\text{a.e.}} \n x$ for some (hence any) $2 \leq n \leq \aleph_0$.
\end{enumerate}
\end{theorem}

\begin{proof}
The equivalence of conditions (1)-(4) is proved as in Theorem \ref{T:appvoic}, with Proposition \ref{T:distance} used in place of Theorem \ref{T:ncomp}.

To see (4) $\Rightarrow$ (5), let $2 \leq n \leq \aleph_0$ and compute
\begin{equation} \label{E:multapp}
x \sim_{\text{a.e.}} \ca y \quad \Rightarrow \quad \n x \sim_{\text{a.e.}} \n \ca y \cong \ca y \sim_{\text{a.e.}} x.
\end{equation}

In \cite[Proof of Corollary 4.3]{H1978b}, Hadwin mentions that for $n=2$, the implication (5) $\Rightarrow$ (2) is a consequence of Voiculescu's theorem.  For the reader's convenience we explicitly prove (5) $\Rightarrow$ (3).  Seeking a contradiction, suppose that $x \sim_{\text{a.e.}} \n x$ for some $2 \leq n \leq \aleph_0$, and that $\text{id}^1$ is not void.  Then
$$x \sim_{\text{a.e.}} \n x \sim_\text{a.e.} \n (\n x) \dots,$$
so by (4) $\Rightarrow$ (3) of Theorem \ref{T:appvoic}, $\text{id}^1$ is $n^k$-divisible for arbitrarily large integer $k$.  But this is impossible, as the range of $\text{id}^1$ contains nonzero finite rank operators.
\end{proof}

\begin{remark} \label{T:now}
We now give an example which shows that if $x$ is a norm limit of $n$-divisibles, the same need not be true for elements of $W^*(x)$.  This was mentioned after Lemma \ref{T:facts}.

Let $x$ be a diagonal operator on $\ell^2$ whose eigenvalues are simple and dense in $[0,1]$, so that $W^*(x)$ contains a rank one projection $p$.  Since $C^*(x) \cap \kom = \{0\}$ and $C^*(p) \cap \kom \ni p$, it follows from Theorem \ref{T:appvoic} ($2 \leq n < \aleph_0$) or \ref{T:appvoic2} ($n = \aleph_0$) that $x$, but not $p$, is a norm limit of $n$-divisible operators.  (It is not hard to argue this directly.)
\end{remark}

\begin{corollary} \label{T:int}
In $\B(\ell^2)$, we have
\begin{equation} \label{E:int1}
\overline{\{\text{$\aleph_0$-divisible operators}\}} = \bigcap_{n < \aleph_0} \overline{ \{\text{$n$-divisible operators}\}}
\end{equation}
On the other hand
\begin{equation} \label{E:int2}
\{\text{$\aleph_0$-divisible operators}\} \subsetneqq \bigcap_{n < \aleph_0} \{\text{$n$-divisible operators}\},
\end{equation}
although this inclusion is dense.
\end{corollary}

\begin{proof}
Equation \eqref{E:int1} follows from the second conditions in Theorems \ref{T:appvoic} and \ref{T:appvoic2}, plus the fact that no compact operator is $n$-divisible for all finite $n$.  The inequality in \eqref{E:int2} results from considering $x \in \B(\ell^2)$ with $W^*(x)'$ of type $\text{II}_1$, while density follows from \eqref{E:int1}.
\end{proof}

Although we will not need Corollary \ref{T:int} in the sequel, we will use

\begin{proposition} \label{T:union}
In $\B(\ell^2)$, we have
\begin{equation} \label{E:unioncl}
\overline{ \bigcup_{n \leq \aleph_0} \{\text{$n$-divisible operators}\}} = \bigcup_{n \leq \aleph_0} \overline{ \{\text{$n$-divisible operators}\}}.
\end{equation}
\end{proposition}


\begin{proof}
We only need to show the inclusion ``$\subseteq$" of \eqref{E:unioncl}.

Let $x \in \B(\ell^2)$, and let $\text{id}$ be the identity representation of $C^*(x)$.  If $\text{id}^1$ is absent, then $C^*(x) \cap \kom = \{0\}$.  By Theorem \ref{T:appvoic2}, $x$ is a norm limit of $\aleph_0$-divisible operators and so belongs to both sides of \eqref{E:unioncl}.  In the remainder of the proof we assume that $\text{id}^1$ is not absent.  This entails that $C^*(x)$ contains a finite rank projection $q$, say of rank $m$.

Suppose that $x$ belongs to the left-hand side of \eqref{E:unioncl}.  Since closure commutes with finite unions,
\begin{equation} \label{E:clunion}
x \in \overline{\{2\text{-divs}\}} \cup \overline{\{3\text{-divs}\}} \cup \dots \cup  \overline{\{m\text{-divs}\}} \cup \overline{ \bigcup_{n > m} \{n\text{-divs}\}}.
\end{equation}
Seeking a contradiction, assume that $x$ does not belong to the right-hand side of \eqref{E:unioncl}.  Then it would have to belong to the union at the far right of \eqref{E:clunion}, and because $q \in C^*(x)$, $q$ would belong to this union as well (by a variant of Lemma \ref{T:facts}(5)).  In particular there must be an operator $d$, $n$-divisible for some $n > m$, with $\|d - q\| = \delta < \frac12$.  By considering the real part, we may assume that $d$ is self-adjoint.

From elementary invertibility considerations $\text{sp}(d) \subseteq [-\delta, \delta] \cup [1-\delta, 1+ \delta]$.  Setting $p = \chi_{[1-\delta, 1+\delta]}(d)$, we compute
\begin{equation} \label{E:qp}
\|p - qp\| \leq \|p - q\| \|p\| = \|p - q\| \leq \|p - d\| + \|d - q\| \leq 2 \delta < 1.
\end{equation}
Now $\|p - q\| < 1$ implies that $p \neq 0$.  Together with the $n$-divisibility of $p$ (because $p \in C^*(d)$) and $n > m$, this gives
$$\text{rank}(qp) \leq \text{rank}(q) = m < \text{rank}(p),$$
so that $qp$ must have nontrivial kernel in $p\ell^2$.  Thus $\|p - qp\| = 1$, contradicting \eqref{E:qp}.

In passing we note that the distance from $q$ to the $n$-divisible operators ($n > m$) is exactly $\frac12$; consider the $n$-divisible operator $\frac12 I$.
\end{proof}



\begin{remark}
It is typically nontrivial to calculate the exact distance from a given $x$ to the $n$-divisible operators.  For $x$ compact and $n= \aleph_0$, we solved this in Proposition \ref{T:distance}.  For $x$ self-adjoint, an answer can be deduced from the main results of \cite{AD}, but we have no need to present such an expression here.

In the proof of Proposition \ref{T:union} we determined that $x$ could not belong to the closure of the $n$-divisible operators, but we did not obtain any estimate of the distance.  Lower bounds are available, at least in theory, by using the ``noncommutative continuous functions" introduced by Hadwin.  These are appropriate limits of noncommutative polynomials; see \cite{H1978a,HKM}.  Revisiting the proof in this light, we have that $q \in C^*(x)$ implies $q = \varphi(x)$ for some noncommutative continuous function $\varphi$.  Continuity means that there is $\delta > 0$ such that for all $y \in \B(\ell^2)$,
$$\|y - x\| < \delta \quad \Rightarrow \quad \|\varphi(y) - q\| = \|\varphi(y) - \varphi(x)\| < \frac12.$$
Still assuming $n > m$, for $\|y - x\| < \delta$ we have
\begin{align*}
\text{dist}(\varphi(y), n\text{-divs}) &\geq \text{dist}(q, n\text{-divs}) - \|\varphi(y) - q\| \\ &= \frac12 - \|\varphi(y) - q\| \\ &> 0,
\end{align*}
which implies as before that $y$ cannot be a norm limit of $n$-divisibles.  Thus the distance from $x$ to the $n$-divisibles is at least $\delta$.

The preceding paragraph bears some resemblance to the proof of \cite[Theorem 2.10]{H1977}.
\end{remark}

Putting Proposition \ref{T:union} together with the implication (1) $\to$ (4) in Theorems \ref{T:appvoic} and \ref{T:appvoic2}, we obtain

\begin{theorem} \label{T:ue}
If an operator $x \in \B(\ell^2)$ is a norm limit of divisible operators, then it is a norm limit of unitary conjugates of a single divisible operator.
\end{theorem}

The skeptical reader may wonder if this is part of a larger and simpler truth, namely that norm limits of unitarily invariant sets in $\B(\ell^2)$ must be approximately equivalent to a member of the set.  A counterexample is given by the irreducible operators, which are norm dense (\cite{Hal1968} or the one-page paper \cite{RR}).  It is easy to check that no irreducible operator is approximately equivalent to a rank one projection.  (Voiculescu's theorem implies that if $x \in \B(\ell^2)$ is approximately equivalent to an irreducible operator, $C^*(x) \cap \kom$ is either $\{0\}$ or $\kom$.)

\smallskip

Theorem \ref{T:typeI}(2) already ruled out the norm density of the $n$-divisible operators in $\B(\ell^2)$, for any $n$.  After a lemma, we will establish a stronger result.

\begin{lemma} \label{T:cont1}
For any $x \in \B(\ell^2)$ and $\eps > 0$, there is $y$ such that $\|y - x\| < \eps$ and $C^*(y)$ contains a rank one projection.
\end{lemma}

\begin{proof}
It goes back to Weyl (\cite{We}) that any self-adjoint operator can be perturbed by an arbitrarily small self-adjoint compact operator to become diagonal.  Apply this to $x_1$, the real part of $x$, finding $k$ such that $\|k\| < \frac{\eps}{2}$ and $x_1 + k$ is diagonal.

Now choose any eigenvalue $\lambda$ for $x_1 + k$ and let $p$ be a rank one projection under $\chi_{\{\lambda\}}(x_1 + k)$.  The operator
$$y = x + k + \frac{\eps}{2} \left[\chi_{(\lambda - (\eps/4), \lambda + (\eps/4))}(x_1 + k) - p \right]$$
has the property that $\chi_{(\lambda - (\eps/4), \lambda + (\eps/4))}(\text{Re }y) = p$, so that $p \in C^*(y)$.  Furthermore
$$\|y - x \| \leq \|k\| + \left\|\frac{\eps}{2} \left[\chi_{(\lambda - (\eps/4), \lambda + (\eps/4))}(x_1 + k) - p \right] \right\| \leq \frac{\eps}{2} + \frac{\eps}{2} = \eps. \qedhere$$
\end{proof}

\begin{theorem} \label{T:nowhere}
The set of divisible operators in $\B(\ell^2)$ is nowhere dense in the norm topology.
\end{theorem}

\begin{proof}
By Lemma \ref{T:cont1} any open ball contains an element $y$ such that $C^*(y)$ contains a rank one projection.  According to Theorems \ref{T:appvoic} and \ref{T:appvoic2}, $y$ is not in the closure of the $n$-divisible operators for any $n$.  By Proposition \ref{T:union}, $y$ is not in the closure of all the divisible operators.
\end{proof}

\smallskip

What about norm density of the $n$-divisible operators in von Neumann algebras of types II and III?  The results of Section \ref{S:op} show that in some $\text{II}_1$ factors, the $n$-divisible operators are not even $\sigma$-strongly dense, but this is all we can say at this point.  It would be interesting to decide the norm density of the $n$-divisible operators in $\R$.  Similarly to Remark \ref{T:mart}, one can apply a martingale theorem (\cite[Theorem 8]{D-N}) to any McDuff factor $(\M, \tau)$ and conclude that any operator in $\M$ is the \textit{almost uniform} limit of $n$-divisible operators.  (This means that for any $\eps > 0$, there is a projection $p \in \M$ and $n$-divisible operators $\{x_n\}$ with $\tau(p) < \eps$ and $\|(x - x_n)p\| \to 0$.)

For any von Neumann algebra, one may measure the size of the norm closure by an invariant analogous to \eqref{E:mcdindex}: just replace $\|\cdot\|_2$ with the operator norm.  All the previous comments (including the author's ignorance) apply to this variation.

Many approximation problems from operator theory are unexplored in the larger context of von Neumann algebras.  Techniques and answers may lend insight into the local structure of the algebras themselves, as in Theorem \ref{T:mcduff}, or even provide useful invariants.  Here is a basic example related both to this paper and to von Neumann algebraic analogues of Voiculescu's theorem.  Say that $x \in \M$ is \textit{reducible} if $W^*(x)' \cap \M \neq \C$ -- \textit{are the reducible operators norm dense in a factor of type II or III?}


\section{Convexity of $\sigma$-weakly closed unitary orbits} \label{S:convex}

In this section we study the possible convexity of $\worb$, where $x \in \M$.  (We remind the reader that the $\sigma$-weak topology is the weak* topology on $\M$.  Actually the Banach space weak topology is also covered by the results below; see Remark \ref{T:bw}.)  For $\M=\B(\ell^2)$, there are some descriptions of $\worb$ in the literature, probably the most notable being Hadwin's characterization as the set of \textit{approximate compressions} of $x$ (\cite[Theorem 4.4(3)]{H1987}):
\begin{equation} \label{E:acomp}
\worb = \overline{\{v^* x v \mid v^* v = 1\}}^{\|\|}.
\end{equation}
It also follows from \cite[Proposition 3.1(3) and Theorem 2.4]{H1987} and \cite[Theorem 2(1)]{HL} that
\begin{align}
\label{E:crn} \worb = \{\varphi(&x) \mid \: \varphi : C^*(x) \to \B(\ell^2) \text{ unital, completely} \\ \notag &\text{positive, and completely rank-nonincreasing}\}.
\end{align}
(The map $\varphi$ is \textit{completely rank-nonincreasing} if
$$\text{id}_n \otimes \varphi: \Mn \otimes C^*(x) \to \Mn \otimes \B(\h)$$
is rank-nonincreasing for all finite $n$.)  In another direction, Kutkut (\cite[Theorem 1.1]{Ku1}) showed that if $x$ is a contraction whose spectrum contains the unit circle, then $\worb$ is the closed unit ball of $\B(\ell^2)$.  He later extended this to certain operators with convex spectral sets (\cite{Ku2}).  Note that the closed unit disk is a spectral set for any contraction, by von Neumann's inequality.

In general von Neumann algebras most of the attention has focused on the closed convex hull $\overline{\text{conv} (\U(x))}^{\|\|}$.  From among the substantial literature, we only mention two results here.  Dixmier's averaging theorem (\cite{Di1949}) establishes that $\overline{\text{conv} (\U(x))}^{\|\|}$ always intersects the center of $\M$.  And assuming that $x$ is self-adjoint and $\M$ has separable predual, Hiai and Nakamura characterized $\overline{\text{conv} (\U(x))}^{\|\|}$ spectrally and proved that it equals $\overline{\text{conv} (\U(x))}^{\sigma-w}$ (\cite{HN}).  So in some cases where we can verify the convexity of $\worb$, we may actually deduce
$$\worb = \overline{\text{conv} (\U(x))}^{\sigma-w} = \overline{\text{conv} (\U(x))}^{\|\|}.$$
This means, for example, that one can do ``Dixmier averaging" without any averaging\dots ${}$ if one is content to approximate in the $\sigma$-weak topology.  One might compare this with \cite[Corollary 6.6]{S2005}, where it is shown that $\overline{\U(x)}^{\sigma-s}$ (but not necessarily $\overline{\U(x)}^{\sigma-s^*}$) intersects the center whenever $\M$ is properly infinite.

Outside of $\B(\h)$ the only descriptions we know of $\worb$ were obtained in recent work with Akemann (\cite{AS}), and they apply exclusively to self-adjoint $x$.  They do show that appropriate generalizations of \eqref{E:acomp} and \eqref{E:crn}, replacing ``rank" by the equivalence class of the range projection, do \textit{not} remain valid.  Concerning convexity, they give

\begin{theorem} \label{T:as} $($\cite{AS}$)$
Let $x$ be a self-adjoint element of a factor $\M$.  If $\M$ is type II or III, then $\worb$ is convex.  If $\M$ is type I, $\worb$ is convex if and only if the spectrum and essential spectrum of $x$ have the same minimum and maximum.
\end{theorem}

At present our only examples of nonconvex $\worb$ are in factors of type I.  We would be interested to know if this can happen in other factors.

The main goal of the section is to prove

\begin{theorem} \label{T:conv}
Let $x$ belong to the norm closure of the $\aleph_0$-divisible operators in a von Neumann algebra $\M$.  Then $\worb$ is convex.
\end{theorem}

As mentioned in the Introduction, this result motivated our entire study of divisible operators.  Its converse is not true: there are operators $x$ which are not norm limits of $\aleph_0$-divisible operators, yet $\worb$ is convex.  (Use Theorems \ref{T:appvoic2} and \ref{T:as}.)  And the implication also fails when the norm topology is replaced by an operator topology.  (Use Proposition \ref{T:spinf} and Theorem \ref{T:as}.)  Actually Theorem \ref{T:conv} is somewhat isolated, but it would have a very nice consequence if one could also show that the $\aleph_0$-divisible operators are norm dense in a type III factor.


\begin{lemma} \label{T:isom}
Let $v$ and $x$ belong to a properly infinite von Neumann algebra $\M$, with $v$ an isometry.  Then $v^* x v \in \worb$.
\end{lemma}

\begin{proof}
Let $\{\varphi_j\} \subset \M^+_*$ be a finite subset.  By repeatedly halving the identity, one can find a decreasing sequence of projections $\{p_n\}$ with
$$p_n \sim 1, \; \forall n \qquad \text{and} \qquad \varphi_j(p_n) \leq \frac1n, \; \forall j,n.$$

For each $n$, $v p_n^\perp$ is a partial isometry with right support $p_n^\perp$ and left support $v p_n^\perp v^*$.  Note that
$$1- v p_n^\perp v^* = 1 - vv^* + v p_n v^* \geq v p_n v^* \sim p_n \sim 1,$$
so that $1- v p_n^\perp v^* \sim p_n$.  Letting $w_n$ be a partial isometry with right support $p_n$ and left support $1 - v p_n^\perp v^*$, define $u_n$ to be the unitary operator $v p_n^\perp + w_n$.  Thus
$$u_n - v = w_n - v p_n.$$

Now we use Cauchy-Schwarz to calculate, for any $j$,
\begin{align*}
|\varphi_j(u_n^* x u_n - v^* x v)| &= |\varphi_j(u_n^* x u_n - v^* x u_n) + \varphi_j(v^* x u_n - v^* x v)| \\ &\leq |\varphi_j((w_n - vp_n)^* x u_n)| +  |\varphi_j(v^* x (w_n - vp_n))| \\ &\leq  \varphi_j (p_n (w_n - v)^*(w_n - v) p_n)^{1/2} \varphi_j(u_n^* x^* x u_n)^{1/2} \\ & \quad + \varphi_j(v^*x x^* v)^{1/2} \varphi_j (p_n (w_n - v)^*(w_n - v) p_n)^{1/2} \\ &\leq\varphi_j (p_n (w_n - v)^*(w_n - v) p_n)^{1/2} (2 \|\varphi_j\|^{1/2} \|x\|) \\ &\leq \varphi_j (4 p_n)^{1/2} (2 \|\varphi_j\|^{1/2} \|x\|) \overset{n \to \infty}{\longrightarrow} 0. \qedhere
\end{align*}
\end{proof}

\begin{proof}[Proof of Theorem \ref{T:conv}]
We first show
\begin{equation} \label{E:mean}
\frac{u_1 x u_1^* + u_2 x u_2^*}{2} \in \worb, \qquad \forall u_1, u_2 \in \U(\M).
\end{equation}

Start by assuming that $x$ is $\aleph_0$-divisible.  Then $W^*(x)' \cap \M$ contains $\B(\ell^2)$, so it contains two isometries $v,w$ satisfying $vv^* + ww^* = 1$.  This implies all of the following:
\begin{equation} \label{E:isom}
v^*w = w^*v = 0, \quad v^*xw = w^*x v =0, \quad v^* x v = w^* x w = x.
\end{equation}

Set
\begin{equation} \label{E:defrs}
r = \frac{vu_1 v^* + v u_2 w^*}{\sqrt{2}}, \quad s = \frac{vu_1 v^* - v u_2 w^*}{\sqrt{2}}.
\end{equation}
By computations using \eqref{E:isom}, one verifies that
$$rr^* = s s^*= vv^*,$$
so $r$ and $s$ are partial isometries, and moreover that
$$r^*r + s^*s = 1.$$
The complements of the left and right supports of $r$ are equivalent:
$$1-rr^* = 1-vv^* = ww^* \sim w^*w = 1 = v^*v \sim vv^* = ss^* \sim s^*s = 1 - r^*r.$$
This means that $r$ can be extended to a unitary $y$, i.e.
\begin{equation} \label{E:defy}
r = r r^* y = v v^* y.
\end{equation}

Using Lemma \ref{T:isom}, we compute
\begin{align*}
\worb &= \overline{\U(yxy^*)}^{\sigma-w} \ni v^* y x y^* v \\ &= v^* v v^* y x y^* v v^* v \\ &\overset{\eqref{E:defy}}{=} v^* r x r^* v \\ &\overset{\eqref{E:defrs}}{=} \frac{u_1 v^* x v u_1^* +  u_1 v^* x w u_2^* + u_2 w^* x v u_1^* + u_2 w^* x w u_2^*}{2} \\ &\overset{\eqref{E:isom}}{=} \frac{u_1 x u_1^* + u_2 x u_2^*}{2}.
\end{align*}

Now we suppose $x$ to be a norm limit of $\aleph_0$-divisibles $\{x_n\}$, as in the statement of the theorem.  For unitaries $u_1, u_2, y$ and a finite set $\{\varphi_j\} \subset \M^+_*$,
\begin{align*}
\left|\varphi_j\left(\frac{u_1 x u_1^* + u_2 x u_2^*}{2} - y xy^*\right)\right| \leq &\left|\varphi_j\left(\frac{u_1 (x-x_n) u_1^* + u_2 (x - x_n) u_2^*}{2}\right)\right| \\ &+  \left|\varphi_j\left(\frac{u_1 x_n u_1^* + u_2 x_n u_2^*}{2} - yx_n y^*\right)\right| \\ &+ |\varphi_j( y (x_n - x) y^*)|.
\end{align*}
We can guarantee that this is small for all $\varphi_j$ by first choosing $n$ to bound the first and third terms, then choosing $y$ as in the first part of the proof to bound the second.  This establishes \eqref{E:mean}.

From \eqref{E:mean} it follows that $\worb \supseteq \text{conv} (\U(x))$.  Then
$$\worb \supseteq \overline{\text{conv} (\U(x))}^{\sigma-w} \supseteq \worb,$$
implying equality.  It is an easy general fact that the closure of a convex set is convex, as long as the map $(\xi,\eta) \mapsto \frac{\xi + \eta}{2}$ is continuous in the relevant topology.  Thus $\worb$ is convex, finishing the proof.
\end{proof}

\begin{remark} \label{T:bw}
The preceding lemma and theorem are also true for the Banach space weak topology ($\sigma(\M, \M^*)$-topology); just choose the set $\{\varphi_j\}$ from $\M_+^*$.

\end{remark}

\appendix

\section{Quotients of operators by cardinals}

We first explain what is meant here by ``dividing an operator by a cardinal."  Given $x \in \M$, a quotient by $n$ is $y \in \N$ satisfying
$$\n (y \in \N) \cong (x \in \M).$$
The existence of a solution is equivalent to the $n$-divisibility of $x$.  As mentioned in the Introduction, uniqueness only becomes meaningful once we agree to identify isomorphic operators, as follows.  Note that $\cong$ is an equivalence relation on operators in von Neumann algebras, and write equivalence classes with brackets, e.g. $[x \in \M]$.  Since amplification is well-defined on equivalence classes (Lemma \ref{T:facts}(1)), we may also consider the equation
\begin{equation} \label{E:abbr}
\n [y \in \N] = [x \in \M].
\end{equation}
For $n$ finite, the solution to \eqref{E:abbr} is always unique (if it exists).  The main goals of this appendix are to explain why this is true and to discuss several variations of interest.

In everything that follows, operators may be replaced with *-homo\-mor\-phisms of $C^*$-algebras.

\subsection{Initial comments} ${}$

\textbf{1.} The first issue in \eqref{E:abbr} is really to identify the algebra $\N$ (up to isomorphism).  We write $[\M]$ for the isomorphism class of $\M$, and $[\M]^n$ for $[\Mn \otimes \M]$.  Then the algebra in question is a solution to
\begin{equation} \label{E:abbr2}
[\N]^n = [\M].
\end{equation}
For factors on a separable Hilbert space, the study of equation \eqref{E:abbr2} goes all the way back to Murray and von Neumann, who wrote it as $\overline{\textbf{N}}^p = \overline{\textbf{M}}$.  Their results (\cite[Section 2.6]{MvN4}) are subsumed in Lemma \ref{T:23div} and Proposition \ref{T:quo}(1).

\smallskip

\textbf{2.} Although we think of the maps $[\M] \mapsto [\M]^n$ and $[x] \mapsto \n [x]$ as multiplications, they do not arise by iterating some kind of sum operation.  Indeed, if there were a sum ``+" satisfying $[\M] \text{``+"} [\M] = [\M]^2$, what would $[\M] \text{``+"} [\N]$ be?  In this sense algebras of the form $\B(\h)$ are very special, since one can take $[\B(\h_1)] \text{``+"} [\B(\h_2)] = [\B(\h_1 \oplus \h_2)]$.  We now admit that the opening paragraph of this paper is somewhat disingenuous.

At the level of operators, the situation is even worse.  In general one cannot form the diagonal sum of a pair of classes, \textit{even from the same algebra}: a pair $[x \in \M], [y \in \M]$ does not determine a well-defined class $[(\begin{smallmatrix} x & 0 \\ 0 & y \end{smallmatrix}) \in \mathbb{M}_2 \otimes \M]$.  As we explain elsewhere, this can be attributed to the existence of automorphisms which are not locally inner (\cite[Section 3]{S2006loc}).

Of course the usual direct products are defined on isomorphism classes of algebras and operators, i.e. $[x \in \M] \oplus [y \in \N] = [x \oplus y \in \M \oplus \N]$.  By iteration, they give rise to the multiplications
\begin{align}
\label{E:cmult1} [\M] &\mapsto [\ell^\infty_n \overline{\otimes} \M],\\
\label{E:cmult2} [x \in \M] &\mapsto [1 \otimes x \in \ell^\infty_n \overline{\otimes} \M].
\end{align}
One may further say that $(y \in \N)$ is ``centrally $n$-divisible" if it is isomorphic to an output of \eqref{E:cmult2}.  But this property is not nicely characterized in $W^*(y)$ or its relative commutant, as the intertwiners which indicate multiplicity lie outside of $\N$.  

The substitution of $\Mn$ for $\ell^\infty_n$ suggests that the maps $[\M] \mapsto [\M]^n$ and $[x] \mapsto \n [x]$ should be considered \textit{quantized} multiplications.

\smallskip

\textbf{3.} In his 1955 book \cite{Ka}, Kaplansky posed three ``test problems" for abelian groups and suggested their possible merit for other  mathematical structures with sum and isomorphism.  (Only the first two problems made it into the second edition of the book.)  The second test problem is this: if $a \oplus a \cong b \oplus b$, must $a \cong b$?  Among the substantial literature on these problems, Kadison and Singer answered them affirmatively in the context of unitary equivalence for Hilbert space operators two years later.  Their result says that in type I factors, $\text{\textcircled{\smaller[3] 2 \larger[3]}} [a] = \text{\textcircled{\smaller[3] 2 \larger[3]}} [b]$ implies $[a]=[b]$.  In other words, the map $[x] \mapsto \text{\textcircled{\smaller[3] 2 \larger[3]}} [x]$ has a (partially-defined) inverse.  See Azoff (\cite{Az}) for other operator theoretic results and references concerning the test problems.  In particular \cite{Az} answers the test problems affirmatively for the direct sum of von Neumann algebras, so that \eqref{E:cmult1} also has an inverse when $n=2$.  Of course any finite $n$ is also suitable; the intrepid reader may go on to show the existence of a ``central" operator quotient by inverting \eqref{E:cmult2}.

\subsection{Uniqueness of finite quotients} It suffices to prove

\begin{proposition} \label{T:quo}
Let $n$ be finite, and let $x \in \M$ and $y \in \N$ be elements of von Neumann algebras.
\begin{enumerate}
\item If $\Mn \otimes \M \cong \Mn \otimes \N$, then $\M \cong \N$.
\item If $\n x \cong \n y$, then $x \cong y$.
\end{enumerate}
\end{proposition}

\begin{proof}[Sketch of proof] Statement (1) is a consequence of (2).  As we just mentioned, (2) was proved for $n=2$ and type I factors in \cite[Theorem 1]{KS}.  The essence of the following argument is the same.

Fix
$$\pi: \Mn \otimes \M \overset{\sim}{\to} \Mn \otimes \N, \qquad \pi(1 \otimes x) = 1 \otimes y.$$
Here is a suitable chain of isomorphisms, with explanations afterward:
\begin{align*}
x \in \M &\cong e_{11} \otimes x \in (e_{11} \otimes 1_\M)(\Mn \otimes \M)(e_{11} \otimes 1_\M) \\
&\cong \pi(e_{11} \otimes 1_\M)(1 \otimes y) \in \pi(e_{11} \otimes 1_\M)(\Mn \otimes \N) \pi(e_{11} \otimes 1_\M) \\
&\cong (e_{11} \otimes 1_\N)(1 \otimes y) \in (e_{11} \otimes 1_\N)(\Mn \otimes \N) (e_{11} \otimes 1_\N) \\
&\cong y \in \N.
\end{align*}

The first and fourth isomorphisms are clear.  The second isomorphism is an application of $\pi$, using $e_{11} \otimes x = (e_{11} \otimes 1_\M)(1 \otimes x)$.

For the third, first note that $e_{11} \otimes 1_\M$ commutes with $1 \otimes x$, so $\pi(e_{11} \otimes 1_\M)$ commutes with $\pi(1 \otimes x) = 1 \otimes y.$  Then both $\pi(e_{11} \otimes 1_\M)$ and $e_{11} \otimes 1_\N$ are projections in $W^*(1 \otimes y)' \cap (\Mn \otimes \N)$ which solve the equation
\begin{equation} \label{E:dimeq}
\underbrace{[p] + [p] + \dots + [p]}_{n \text{ times }} = [1]
\end{equation}
in the dimension theory for $W^*(1 \otimes y)' \cap (\Mn \otimes \N)$, and this implies that they are Murray-von Neumann equivalent in $W^*(1 \otimes y)' \cap (\Mn \otimes \N)$.  (The \textit{dimension theory} for a von Neumann algebra $\M$ is the quotient $\pM$, where $\p(\M)$ is the set of the projections in $\M$ and $\sim$ is Murray-von Neumann equivalence.  Among its many features is a partially-defined addition for arbitrarily large sets of summands.  See \cite[Section 2]{S2006dim} for an overview.)  The third isomorphism can then be had by conjugating by a partial isometry in $W^*(1 \otimes y)' \cap (\Mn \otimes \N)$ which goes from $\pi(e_{11} \otimes 1_\M)$ to $e_{11} \otimes 1_\N$.
\end{proof}

Although it looks innocuous, Proposition \ref{T:quo}(1) does not hold for $C^*$-algebras!  The first example was given in \cite{Pl}, and \cite{Ko} contains a more systematic study.

Proposition \ref{T:quo} also fails for infinite $n$.  For example, note that a projection in $\B(\ell^2)$ with infinite rank and corank is an $\aleph_0$-multiple of any nontrivial projection on a separable (possibly finite-dimensional) Hilbert space.  Just as for cardinals, division by an infinite quantity is problematic.  We do, however, have the implications
$$[\N]^n = [\M] \Rightarrow [\M]^n = [\M] \quad \text{and} \quad \n [y] = [x] \Rightarrow \n [x] = [x].$$
Their proofs are alike; for the second, assume $\n y \cong x$ and compute
$$\n x \cong \n (\n x) \cong \n y \cong x.$$
This means that for $n$ infinite and $x$ $n$-divisible, the equation $\n [y] = [x]$ always has the solution $[y] = [x]$.  Typically there are other solutions, but not always (for instance, $n = \aleph_0$ and $x$ the identity of a $\sigma$-finite type III factor).

The property
\begin{equation} \label{E:ss}
[x] = \n [x]
\end{equation}
may be thought of as a ``self-similarity."  For $n$ infinite, \eqref{E:ss} is no stronger than $n$-divisibility, as we just mentioned.   For $n$ finite, by repeated substitution \eqref{E:ss} entails that $x$ is $n^k$-divisible for any natural $k$, or equivalently, that $W^*(x)' \cap \M$ lacks a finite type I summand.  But the converse to this implication does not generally hold.  For example, the identity of a $\text{II}_1$ factor $\M$ satisfies \eqref{E:ss} if and only if $\M \cong \Mn \otimes \M$, which is not always true.  


\subsection{Generalized amplifications of operators} Readers familiar with von Neumann algebras will not be surprised to hear that for some $[x \in \M]$, the map $[x] \mapsto \n [x]$ makes sense for non-integer values of $n$.  So for example one may sometimes amplify $[x \in \M]$ by $\sqrt{2}$, thinking of this as the isomorphism class of the quantum direct sum of $\sqrt{2}$ copies of $x$.

In fact, in the broadest context, the parameter may be chosen from the Murray-von Neumann equivalence classes of projections in amplifications of $W^*(x)' \cap \M$.  For a projection $p \in \mathbb{M}_k \otimes (W^*(x)' \cap \M) \subseteq \mathbb{M}_k \otimes \M$, we define the \textit{(generalized) amplification of $[x \in \M]$ by the dimension $[p]$} to be
$$\text{\textcircled{\smaller[5]$[p]$\larger[5]}} [x \in \M] \overset{\text{def}}{=} [p(1 \otimes x) \in p(\mathbb{M}_k \otimes \M)p].$$
The parameter $[p]$ looks dishearteningly non-numerical, but by dimension theory it can be identified with a cardinal-valued function on the spectrum of the center of $W^*(x)' \cap \M$ (\cite{To}).  (On the spectrum of the type II summand, the function may also take values in the positive reals.)

This allows us to unify division and multiplication, as we now illustrate with a simple example.  Consider \eqref{E:abbr} under the assumption that $W^*(x)' \cap \M$ and $n$ are finite, with the identity of $W^*(x)' \cap \M$ $n$-divisible.  Let $[p] \in (\p(W^*(x)' \cap \M)/\sim)$ be the unique solution to \eqref{E:dimeq}.  Then $[p]$ is characterized as the set of projections in $W^*(x)' \cap \M$ whose image under the canonical dimension function (or center-valued trace) for $W^*(x)' \cap \M$ is exactly $\frac1n$.  On its domain the operation $\text{\textcircled{\smaller[5]$[p]$\larger[5]}}$ is inverse to $\n$, so one considers it as ``division by $n$," solving \eqref{E:abbr} for $[y]$ by applying it to both sides.

It is actually not too much trouble to set up an algebraic calculus for amplifications in which only dimensions are used.  But in general the incorporation of cardinals requires some unwieldy extra bookkeeping, because dimension functions are not unital in infinite algebras (so that the map from cardinals to dimensions is many-to-one), and are not even canonical in $\text{II}_\infty$ algebras.  We do not give the details here, but we point out that Ernest worked out a version of this theory for $\B(\ell^2)$ (\cite[Chapter 4]{E}).  He used no cardinals higher than $\aleph_0$, and he only considered dimensions with full central support.  (Modulo the cardinality restriction, these correspond to the \textit{coupling functions} for $W^*(x)' \cap \M$.)  This produces a useful subset of the amplifications of $x \in \B(\ell^2)$:
\begin{equation} \label{E:qe}
\{y \in \B(\ell^2) \mid \ca y \cong \ca x\}.
\end{equation}

Ernest called the relation $\ca y \cong \ca x$ \textit{quasi-equivalence}, so that \eqref{E:qe} is the quasi-equivalence class of $x$.  This relation might also be considered ``stable equivalence" or even a sort of Morita equivalence for operators, as coupling functions are invariants of representations.  In $\B(\ell^2)$ or other $\sigma$-finite von Neumann algebras we prefer to call \eqref{E:qe} the \textit{genus} of $x$, following terminology set up long ago by Murray and von Neumann for analogous equivalence classes of factors (\cite[Chapter III]{MvN4}).  So once again, they have the last word.



\begin{thebibliography}{HKM}

\bibitem[AA]{AA}
C. A. Akemann and J. Anderson,
\emph{Lyapunov theorems for operator algebras},
Mem. Amer. Math. Soc. \textbf{94} (1991), no. 458.

\bibitem[AS]{AS}
C. A. Akemann and D. Sherman,
\emph{The weak*-closed unitary orbit of a self-adjoint operator in a factor},
in preparation.

\bibitem[A1]{A1976}
W. Arveson,
An Invitation to $C^*$-Algebras,
Graduate Texts in Mathematics, vol. 39, Springer-Verlag, New York, 1976.

\bibitem[A2]{A1977}
W. Arveson,
\emph{Notes on extensions of $C^*$-algebras},
Duke Math. J. \textbf{44} (1977), 329--355.

\bibitem[Az]{Az}
E. Azoff,
\emph{Test problems for operator algebras},
Trans. Amer. Math. Soc. \textbf{347} (1995), 2989--3001.

\bibitem[AD]{AD}
E. Azoff and C. Davis,
\emph{On distances between unitary orbits of selfadjoint operators},
Acta Sci. Math. (Szeged) \textbf{47} (1984), 419--439.

\bibitem[Bl]{Bl}
B. Blackadar,
Operator Algebras: Theory of $C^*$-algebras and von Neumann Algebras,
Springer-Verlag, Berlin, 2006.

\bibitem[B]{B}
D. Bisch,
\emph{On the existence of central sequences in subfactors},
Trans. Amer. Math. Soc. \textbf{321} (1990), 117--128.


\bibitem[BKR]{BKR}
B. Blackadar, A. Kumjian, and M. R\o rdam,
\emph{Approximately central matrix units and the structure of noncommutative tori},
$K$-Theory \textbf{6} (1992), 267--284.

\bibitem[Bu]{Bu}
D. Bures,
\emph{Abelian subalgebras of von Neumann algebras},
Mem. Amer. Math. Soc. (1971), no. 110.


\bibitem[C1]{C1976}
A. Connes,
\emph{Classification of injective factors: cases $\text{II}_1$, $\text{II}_\infty$, $\text{III}_\lambda$, $\lambda \neq 1$},
Ann. of Math. (2) \textbf{104} (1976), 73--115.

\bibitem[C2]{C1980}
A. Connes,
\emph{A factor of type $\text{II}_1$ with countable fundamental group},
J. Operator Theory \textbf{4} (1980), 151--153.

\bibitem[CJ]{CJ}
A. Connes and V. Jones,
\emph{Property $T$ for von Neumann algebras},
Bull. London Math. Soc. \textbf{17} (1985), 57--62.

\bibitem[D-N]{D-N}
N. Dang-Ngoc,
\emph{Pointwise convergence of martingales in von Neumann algebras},
Israel J. Math. \textbf{34} (1979), 273--280.


\bibitem[DH]{DH}
H. Ding and D. Hadwin,
\emph{Approximate equivalence in von Neumann algebras},
Sci. China Ser. A \textbf{48} (2005), 239--247.

\bibitem[D]{Di1949}
J. Dixmier,
\emph{Les anneaux d'op\'{e}rateurs de classe finie},
Ann. Sci. \'{E}cole Norm. Sup. (3) \textbf{66} (1949), 209--261.

\bibitem[DL]{DL}
J. Dixmier and E. C. Lance,
\emph{Deux nouveaux facteurs de type $\text{II}_1$},
Invent. Math. \textbf{7} (1969), 226--234.

\bibitem[Dy]{Dy}
H. A. Dye,
\emph{The unitary structure in finite rings of operators},
Duke Math. J. \textbf{20} (1953), 55--69.

\bibitem[Ef]{Ef}
E. Effros,
\emph{The Borel space of von Neumann algebras on a separable Hilbert space},
Pacific J. Math. \textbf{15} (1965), 1153--1164.

\bibitem[El]{El}
G. A. Elliott,
\emph{On derivations of $AW^*$-algebras},
T\^{o}hoku Math. J. (2) \textbf{30} (1978), 263--276.

\bibitem[E]{E}
J. Ernest,
\emph{Charting the operator terrain},
Mem. Amer. Math. Soc. \textbf{6} (1976), no. 171.

\bibitem[FGL]{FGL}
J. Fang, L. Ge, and W. Li,
\emph{Central sequence algebras of von Neumann algebras},
Taiwanese J. Math. \textbf{10} (2006), 187--200.

\bibitem[Ga]{Ga}
P. Gajendragadkar,
\emph{Norm of a derivation on a von Neumann algebra},
Trans. Amer. Math. Soc. \textbf{170} (1972), 165--170.


\bibitem[GH]{GH}
L. Ge and D. Hadwin,
\emph{Ultraproducts of C*-algebras},
in: Recent Advances in Operator Theory and Related Topics (Szeged, 1999), pp. 305--326,
Oper. Theory Adv. Appl., vol. 127, Birkh\"{a}user, Basel, 2001.

\bibitem[GP]{GP}
L. Ge and S. Popa,
\emph{On some decomposition properties for factors of type $\text{II}_1$},
Duke Math. J. \textbf{94} (1998), 79--101.

\bibitem[GS]{GS}
L. Ge and J. Shen,
\emph{Generator problem for certain property $T$ factors},
Proc. Natl. Acad. Sci. USA \textbf{99} (2002), 565--567.

\bibitem[G]{G}
T. Giordano,
personal communication.

\bibitem[GrP]{GrP}
K. Grove and G. K. Pedersen,
\emph{Diagonalizing matrices over $C(X)$},
J. Funct. Anal. \textbf{59} (1984), 65--89.


\bibitem[HW]{HW1}
U. Haagerup and C. Winsl\o w,
\emph{The Effros-Mar\'{e}chal topology in the space of von Neumann algebras},
Amer. J. Math. \textbf{120} (1998), 567--617.



\bibitem[H1]{H1977}
D. Hadwin,
\emph{An operator-valued spectrum},
Indiana Univ. Math. J. \textbf{26} (1977), 329--340.

\bibitem[H2]{H1978a}
D. Hadwin,
\emph{Continuous functions of operators: a functional calculus},
Indiana Univ. Math. J. \textbf{27} (1978), 113--125.

\bibitem[H3]{H1978b}
D. Hadwin,
\emph{An asymptotic double commutant theorem for $C^*$-algebras},
Trans. Amer. Math. Soc. \textbf{244} (1978), 273--297.


\bibitem[H4]{H1987}
D. Hadwin,
\emph{Completely positive maps and approximate equivalence},
Indiana Univ. Math. J. \textbf{36} (1987), 211--228.

\bibitem[H5]{H1998}
D. Hadwin,
\emph{Free entropy and approximate equivalence in von Neumann algebras}, in: Operator Algebras and Operator Theory (Shanghai, 1997), pp. 111-131, Contemp. Math., vol. 228, Amer. Math. Soc., Providence, 1998.

\bibitem[HKM]{HKM}
D. Hadwin, L. Kaonga, and B. Mathes,
\emph{Noncommutative continuous functions},
J. Korean Math. Soc. \textbf{40} (2003), 789--830.

\bibitem[HL]{HL}
D. Hadwin and D. Larson,
\emph{Completely rank-nonincreasing linear maps},
J. Funct. Anal. \textbf{199} (2003), 210--227.

\bibitem[Ha1]{Hal1951}
P. R. Halmos,
Introduction to Hilbert Space and the Theory of Spectral Multiplicity,
Chelsea, New York, 1951.

\bibitem[Ha2]{Hal1968}
P. R. Halmos,
\emph{Irreducible operators},
Michigan Math. J. \textbf{15} (1968) 215--223.

\bibitem[Ha3]{Hal1970}
P. R. Halmos,
\emph{Ten problems in Hilbert space},
Bull. Amer. Math. Soc. \textbf{76} (1970), 887--933.


\bibitem[Hal]{Halp1977}
H. Halpern,
\emph{Essential central range and selfadjoint commutators in properly infinite von Neumann algebras},
Trans. Amer. Math. Soc. \textbf{228} (1977), 117--146.

\bibitem[HN]{HN}
F. Hiai and Y. Nakamura,
\emph{Closed convex hulls of unitary orbits in von Neumann algebras},
Trans. Amer. Math. Soc. \textbf{323} (1991), 1--38.

\bibitem[K1]{K1957}
R. V. Kadison,
\emph{Unitary invariants for representations of operator algebras},
Ann. of Math. (2) \textbf{66} (1957), 304--379.

\bibitem[K2]{K1984}
R. V. Kadison,
\emph{Diagonalizing matrices},
Amer. J. Math. \textbf{106} (1984), 1451--1468.

\bibitem[KK]{KK}
R. V. Kadison and D. Kastler,
\emph{Perturbations of von Neumann algebras I: stability of type},
Amer. J. Math. \textbf{94} (1972), 38--54.

\bibitem[KR]{KR}
R. V. Kadison and J. R. Ringrose,
Fundamentals of the Theory of Operator Algebras I,
Graduate Studies in Mathematics, vol. 15, Amer. Math. Soc., Providence, 1997.

\bibitem[KS]{KS}
R. V. Kadison and I. M. Singer,
\emph{Three test problems in operator theory},
Pacific J. Math. \textbf{7} (1957), 1101--1106.

\bibitem[Kaf1]{K}
V. Kaftal,
\emph{On the theory of compact operators in von Neumann algebras I},
Indiana Univ. Math. J. \textbf{26} (1977), 447--457.

\bibitem[Kaf2]{Kaf1991}
V. Kaftal,
\emph{Type decomposition for von Neumann algebra embeddings},
J. Funct. Anal. \textbf{98} (1991), 169--193.

\bibitem[Ka1]{Ka1951}
I. Kaplansky,
\emph{A theorem on rings of operators},
Pacific J. Math. \textbf{1} (1951), 227--232.

\bibitem[Ka2]{Ka}
I. Kaplansky,
Infinite Abelian Groups,
The University of Michigan Press, Ann Arbor, 1955.

\bibitem[Ko]{Ko}
K. Kodaka,
\emph{$C^*$-algebras that are isomorphic after tensoring and full projections},
Proc. Edinb. Math. Soc. (2) \textbf{47} (2004), 659--668.

\bibitem[Ku1]{Ku1}
M. Kutkut,
\emph{Weak closure of the unitary orbit of contractions},
Acta Math. Hungar. \textbf{46} (1985), 255--263.

\bibitem[Ku2]{Ku2}
M. Kutkut,
\emph{Weak closure of the unitary orbits of operators having a convex spectral set},
J. Pure Appl. Sci. \textbf{18} (1985), 271--281.


\bibitem[MM]{MM}
L. W. Marcoux and G. J. Murphy,
\emph{Unitarily-invariant linear spaces in $C^*$-algebras},
Proc. Amer. Math. Soc. \textbf{126} (1998), 3597--3605.

\bibitem[M]{M}
K. Matsumoto,
\emph{Strongly stable factors, crossed products and property $\Gamma$},
Tokyo J. Math. \textbf{11} (1988), 247--268.

\bibitem[McD]{McD}
D. McDuff,
\emph{Central sequences and the hyperfinite factor},
Proc. London Math. Soc. (3) \textbf{21} (1970), 443--461.

\bibitem[MvN]{MvN4}
F. J. Murray and J. von Neumann,
\emph{On rings of operators IV},
Ann. of Math. (2) \textbf{44} (1943), 716--808.

\bibitem[vN1]{vN}
J. von Neumann,
\emph{Zur Algebra der Funktionaloperatoren und Theorie der normalen Operatoren},
Math. Ann. \textbf{102} (1929), 370--427.

\bibitem[vN2]{vN1942}
J. von Neumann,
\emph{Approximative properties of matrices of high finite order},
Portugaliae Math. \textbf{3} (1942), 1--62.


\bibitem[Pl]{Pl}
J. Plastiras,
\emph{$C^*$-algebras isomorphic after tensoring},
Proc. Amer. Math. Soc. \textbf{66} (1977), 276--278.


\bibitem[P1]{P1985}
S. Popa,
\emph{Notes on Cartan subalgebras in type $\text{II}_1$ factors},
Math. Scand. \textbf{57} (1985), 171--188.

\bibitem[P2]{P1995}
S. Popa,
\emph{Free-independent sequences in type $\text{II}_1$ factors and related problems}, in: Recent Advances in Operator Algebras (Orl\'{e}ans, 1992),
Ast\'{e}risque no. 232 (1995), 187--202.


\bibitem[RR]{RR}
H. Radjavi and P. Rosenthal,
\emph{The set of irreducible operators is dense},
Proc. Amer. Math. Soc. \textbf{21} (1969), 256.


\bibitem[Sa]{Sa}
S. Sakai,
The Theory of $W^*$-Algebras, lecture notes, Yale University, 1962.


\bibitem[Sh]{Sh}
J. Shen,
\emph{Singly generated $\text{II}_1$ factors},
arXiv preprint math.OA/0511327.

\bibitem[S1]{S2005}
D. Sherman,
\emph{Unitary orbits of normal operators in von Neumann algebras},
J. Reine Angew. Math. \textbf{605} (2007), 95--132.

\bibitem[S2]{S2006dim}
D. Sherman,
\emph{On the dimension theory of von Neumann algebras},
Math. Scand. \textbf{101} (2007), 123--147.

\bibitem[S3]{S2006loc}
D. Sherman,
\emph{Locally inner automorphisms of operator algebras},
arXiv preprint math.OA/0609735.

\bibitem[S4]{S2008}
D. Sherman,
\emph{Remarks on automorphisms of ultrapowers of $\text{II}_1$ factors},
in preparation.


\bibitem[Sho]{Sho}
K. Shoda,
\emph{Einige S\"{a}tze \"{u}ber Matrizen},
Japan J. Math \textbf{13} (1936), 361--365.

\bibitem[T1]{T}
M. Takesaki,
Theory of Operator Algebras I,
Springer-Verlag, Berlin, 1979.


\bibitem[T2]{T3}
M. Takesaki,
Theory of Operator Algebras III,
Springer-Verlag, Berlin, 2002.

\bibitem[To]{To}
J. Tomiyama,
\emph{Generalized dimension function for $W^*$-algebras of infinite type},
T\^{o}hoku Math. J. (2) \textbf{10} (1958), 121--129.

\bibitem[U]{U}
H. Umegaki,
\emph{Conditional expectation in an operator algebra II},
T\^{o}hoku Math. J. (2) \textbf{8} (1956), 86--100.

\bibitem[V]{V}
D. Voiculescu,
\emph{A non-commutative Weyl-von Neumann theorem},
Rev. Roumaine Math. Pures Appl. \textbf{21} (1976), 97--113.

\bibitem[We]{We}
H. Weyl,
\emph{\"{U}ber beschr\"{a}nkte quadritische Formen deren Differenz vollsteig ist},
Rend. Circ. Mat. Palermo \textbf{27} (1909), 373--392.

\bibitem[W]{W}
F. B. Wright,
\emph{A reduction for algebras of finite type},
Ann. of Math. (2) \textbf{60} (1954), 560--570.

\bibitem[Z]{Z}
L. Zsid\'{o},
\emph{The norm of a derivation in a $W^*$-algebra},
Proc. Amer. Math. Soc. \textbf{38} (1973), 147--150.

\end{thebibliography}
\end{document}